\documentclass[graybox]{svmult}

\usepackage{type1cm}        % activate if the above 3 fonts are
\usepackage{makeidx}         % allows index generatio1
\usepackage{graphicx}        % standard LaTeX graphics tool
\usepackage{multicol}        % used for the two-column index
\usepackage[bottom]{footmisc}% places footnotes at page bottom

\usepackage{newtxtext}       % 
\usepackage[varvw]{newtxmath}       % selects Times Roman as basic font
\usepackage{amsmath}
\usepackage{fancyhdr}
\makeindex             % used for the subject index
\let\Bbb\mathbb
\let\goth\mathfrak

\def\og{\leavevmode\raise.3ex\hbox{$\scriptscriptstyle\langle\!\langle\,$}}

\def\fgf{\/\leavevmode\raise.3ex\hbox{$\scriptscriptstyle\,\rangle\!\rangle$}}

\def\fg{\fgf\ }

\newcommand{\carre}{\qed}

\newcommand{\findem}{\ensuremath\blacksquare}

\def\Demd#1|{\parindent=0pt\par{\sl D\'emonstration d#1}\pointir\parindent=20pt}

\let\petcap\sc

\def\finc{\vskip12pt}

\def\Dem{\parindent=0pt\par{\sl D\'emonstration}\pointir\parindent=20pt}

\def\Demo#1|{\parindent=0pt\par{\sl D\'emonstration #1}\pointir\parindent=20pt}

\makeatletter
\newdimen\indentTh\indentTh=0pt
\def\p@int{{\rm .}}
\def\p@intir{\discretionary{\rm .}{}{\rm .\kern.35em---\kern.7em}}
\def\pointir{\afterassignment\pointir@\let\next=}
\def\pointir@{\ifx\next\par\p@int\else\p@intir\fi\next}
\long\def\Thc#1|#2\finc{\Th{}{#1}{\pointir}{#2}}
\long\def\Th#1#2#3#4{\parindent=\indentTh\par\vskip5pt
{#1}{\petcap #2}{\sl #3}\parindent=20pt{\sl #4\par}\vskip 5pt\parindent=20pt}

\newdimen\indentssec\indentssec=20pt
\newdimen\indentrem\indentssec=0pt

\def\Demdsp#1|{\parindent=0pt\par{\sl D\'emonstration d#1.}\parindent=20pt}
\def\oldstyle{}
\def\build#1_#2^#3{\mathrel{\mathop{\kern 0pt#1}\limits_{#2}^{#3}}}

\def\hdfl#1#2{\smash{\mathop{\hbox to 12mm{\rightarrowfill}}
\limits^{\scriptstyle#1}_{\scriptstyle#2}}}

\def\hdhfl#1#2{\smash{\mathop{\hbox to 12mm{\hookrightarrowfill}}
\limits^{\scriptstyle#1}_{\scriptstyle#2}}}

\def\hgfl#1#2{\smash{\mathop{\hbox to 12mm{\leftarrowfill}}
\limits^{\scriptstyle#1}_{\scriptstyle#2}}}

\def\hghfl#1#2{\smash{\mathop{\hbox to 12mm{\hookleftarrowfill}}
\limits^{\scriptstyle#1}_{\scriptstyle#2}}}

\newcount\secno \secno=0
\newcount\ssecno \ssecno=0
\newcount\sssecno \sssecno=0
\newcount\chapno \chapno=0
\newcount\notenumber \notenumber=1
\newcount\exino \exino=0
\newcount\expno \expno=0
\newdimen\indentsec\indentsec=20pt
\newdimen\indentssec\indentssec=20pt
\newdimen\indentsssec\indentsssec=20pt
\newdimen\indentrem\indentssec=0pt

\newdimen\indentTh\indentTh=0pt

\newdimen\indentth\indentssec=0pt
\def\sectiongen#1#2#3{\parindent=\indentsec\par\vskip .3cm
\vskip 0mm plus -20mm minus 1,5mm\penalty-50
{\bf #1}{\bf #2}{#3}\nobreak\parindent=20pt}%
\def\secc#1|{\sectiongen{}{#1}{\pointir}}
\def\nsecc#1|%
{\global\advance\secno by 1\global\ssecno=0\global\sssecno=0
\sectiongen{\the\secno\ }{#1}{\pointir}}
\def\secp#1|{\sectiongen{}{#1}{}\par}
\def\nsecp#1|%
{\global\advance\secno by 1\global\ssecno=0\global\sssecno=0
\sectiongen{\the\secno\ }{#1}{}\par}
\def\ssectiongen#1#2#3#4{\parindent=\indentssec\par\vskip .2cm
\vskip 0mm plus -20mm minus 1,5mm\penalty-50
{\bf #1}{\sl #2}{\sl #3}{#4}\nobreak\medskip\parindent=20pt}%
\def\ssecc#1|#2{\ssectiongen{}{#1}{\pointir}{#2}}
\def\nssecc#1|#2{\global\advance\ssecno by 1\global\sssecno=0
\ssectiongen{\the\secno.\the\ssecno\ }{#1}{\pointir}{#2}}
\def\ssecp#1|{\ssectiongen{}{#1}{}{}\par}
\def\nssecp#1|{\global\advance\ssecno by 1\global\sssecno=0
\ssectiongen{\the\secno.\the\ssecno\ }{#1}{}{}\par}
\def\sssectiongen#1#2#3#4{\parindent=\indentsssec\par\vskip .2cm
\vskip 0mm plus -20mm minus 1,5mm\penalty-50
{\bf #1}{\sl #2}{\sl #3}{#4}\nobreak\medskip\parindent=20pt}%
\def\sssecc#1|#2{\sssectiongen{}{#1}{\pointir}{#2}}
\def\nsssecc#1|#2{\global\advance\sssecno by 1
\ssectiongen{\the\secno.\the\ssecno.\the\sssecno\ }{#1}{\pointir}{#2}}
\def\sssecp#1|{\sssectiongen{}{#1}{}{}\par}
\def\nsssecp#1|{\global\advance\sssecno by 1
\sssectiongen{\the\secno.\the\ssecno.\the\sssecno\ }{#1}{}{}\par}
\long\def\Th#1#2#3#4{\parindent=\indentTh\par\vskip5pt
{#1}{\petcap #2}{\sl #3}\parindent=20pt{\sl #4\par}\vskip 5pt\parindent=20pt}
\long\def\pTh#1#2#3#4{\parindent=\indentth\par\vskip5pt
{#1}{\small \bf #2}{\small \sl #3}\parindent=20pt{\small \sl #4\par}\vskip 5pt\parindent=20pt}
\long\def\remarque#1#2#3#4{\parindent=\indentrem\par\vskip5pt
{#1}{\small \sl #2}{\small\sl #3}\parindent=20pt{\small#4\par}
\vskip 5pt\parindent=20pt}
\long\def\remarques#1#2#3#4{\parindent=\indentrem\par\vskip5pt
{#1}{\small \sl #2}{\small\sl #3}{\small#4\par}
\vskip 5pt\parindent=20pt}

\long\def\remarquesa#1#2#3#4#5{\parindent=\indentrem\par\vskip5pt
{#1}{\small \sl #2}{\small\sl #3}{\small#4}
\parindent=20pt{\small#5\par}
\vskip 5pt\parindent=20pt}

\long\def\Remarque#1#2#3#4{\parindent=\indentrem\par\vskip5pt
{#1}{ \sl #2}{\sl #3}\parindent=20pt{#4\par}
\vskip 5pt\parindent=20pt}
\long\def\remarquesn#1#2#3#4{\parindent=\indentrem\par\vskip5pt
{\small \sl #1}{\small #2}{\small\sl #3}{\small#4\par}
\vskip 5pt\parindent=20pt}

\long\def\Remarquen#1#2#3#4{\parindent=\indentrem\par\vskip5pt
{ \sl #1}{#2}{\sl #3}\parindent=20pt{#4\par}
\vskip 5pt\parindent=20pt}

\long\def\Thc#1|#2\finc{\Th{}{#1}{\pointir}{#2}}
\long\def\Thnc#1|#2|#3\finnc{\Th{#1}{#2}{\pointir}{#3}}
\long\def\Exic#1|#2\finc{{\global\advance\exino by 1}\Remarquen%
{Exercice }{\the\exino}{ #1\pointir}{#2}}
\long\def\Expc#1|#2\finc{{\global\advance\expno by 1}\Remarquen%
{Exemple }{\the\expno}{ #1\pointir}{#2}}
\long\def\exic#1|#2\finc{{\global\advance\exino by 1}\remarquesn%
{\bf Exercice }{\the\exino}{ #1\pointir}{#2}}
\long\def\expc#1|#2\finc{{\global\advance\expno by 1}\remarquesn%
{Exemple }{\the\expno}{ #1\pointir}{#2}}

\long\def\Ec#1\finc{\Th{}{}{}{#1}}
\long\def\thc#1|#2\finc{\pTh{}{#1}{\pointir}{#2}}%
\long\def\thsnc#1\finc{\pTh{}{}{}{#1}}

\long\def\Thp#1|#2\finp{\Th{}{#1}{\par}{#2}}
\long\def\thp#1|#2\finp{\pTh{}{#1}{\par}{#2}}
\long\def\rmc#1|#2\finc{\remarque{}{#1}{\pointir}{#2}}
\long\def\Rmc#1|#2\finc{\Remarque{}{#1}{\pointir}{#2}}

\long\def\rmp#1|#2\finp{\remarque{}{#1}{\par}{#2}}
\long\def\Rmp#1|#2\finp{\Remarque{}{#1}{\par}{#2}}

\long\def\parc#1\finc{\remarque{}{}{}{#1}}
\long\def\parcs#1\fincs{\remarques{}{}{}{#1}}

\long\def\parcsa#1\fins#2\fincsa{\remarquesa{}{}{}{#1}{#2}}

\def\Rm#1|{\parindent=0pt\par\vskip5pt{\sl #1}\pointir\parindent=20pt}
\def\preuved#1|{\parindent=0pt\par{\sl Preuve d#1}\pointir\parindent=20pt}
\def\Demf{\parindent=0pt\par{\sl Fin de la d\'emonstration}\pointir\parindent=20pt}
\def\cotg{\mathop{\rm cotg}\nolimits}
\def\nDefns{\parindent=0pt\par{\rm\bf D\'efinitions}$\,$\parindent=20pt}
\input boxedeps.tex
\SetDVIWindowEPSFSpecial
\SetEPSFDirectory{}
\SetDefaultEPSFScale{250}
\EmulateRokicki

\global\indentsec=0pt\global\ssecno=0\global\secno=0
\global\ssecno=0\global\secno=0
\def\page {\leaders\hbox to 2mm{\hfil.\hfil}\hfill
}
\def\npage {\vfill\eject \global\setcounter{footnote}{0}}
\def\npg {\global\setcounter{footnote}{0}}
\def\page#1#2{\leaders\hbox to 1mm{\hfil.\hfil}\hfill
\rlap{\hbox to 10mm{\hfill #1}{\ \  p.#2}}\par
}
\def\tpage #1{\leaders\hbox to 2mm{\hfil.\hfil}\hfill
\rlap{\hbox to 15mm{\hfill #1}}\par
}
\def\npage {\vfill\eject \global\setcounter{footnote}{0}}

\begin{document}
\fancyhf{}
\pagestyle{fancy}
%\setlength{\headheight}{12.72972pt}
%\addtolength{\topmargin}{-0.72972pt}
%\renewcommand{\headrulewidth}{0.2pt}

\title*{C'est un petit val qui mousse de fonctions}

\author{{\it par} Alexis Marin, illustrations de Greg McShane}
\maketitle
\vskip-40mm
\hfill {\small\sl \`a Jean Cerf pour son nonante-neuvi\`eme anniversaire.\/}

\vskip20mm
\parc
\centerline{\bf R\'esum\'e}
%\og%
Aux {\it $f$
et
$\theta$%
-fibr\'es\/} %\fg
de Marston Morse on substitue ici les
{\it vals de fonctions\/},
ils se modifient par le
\og trivial abaissement  des minimas\fgf,
se caract\'erisent par la m\'ethode du chemin de Moser
et souplement d\'ecrivent
\og stabilisation de fonctions\fg et
\og r\'eduction au cas des courbes\fgf.

Une excursion passant par une \og extension canonique de nappe transverse\fg
conduira la lectrice, sans technicit\'es de plongement de mod\`eles,
aux
\og Lemmes de Morse\fgf~:\hfill\break  
forme canonique,
abaissement de valeur critique et
\'elimination de paire de points critiques.
\finc
\vskip5mm
\parc
\centerline{\bf Abstract}
Easy  lowering local minima,
after introducing
{\it valley functions\/}
on  smooth mani\-folds gives, 
without gluing technicities,
``M. Morse's lemma's''~: canonical form, moving critical values and 
eliminating pair of critical points theorems,
reducing the last  to its one dimensional case.

\finc
\vskip3mm

{\noindent \small 2020 \sl{Mathematics Subject Classification}~:\hfill\break
Primary: 57R35; Secondary: 57R80, 57-03.\hfill\break
}

\vskip15mm
%\midinsert
%\centerline{\BoxedEPSF{val1.eps }}
%\centerline{}
%\endinsert

%\parc
L'introduction en cort\`ege de d\'efinitions et \'enonc\'es au \S1
se compl\`ete par la postroduction%
\footnote{\small
{\it c.a d.\/} une introduction compr\'ehensible par qui,
sans entrer dans la technicit\'e des appendices,
comprend les quatre \'enonc\'es principaux et leur r\^ole
 au \S2 pour l'\'elimination.\hfill\break
On peut alors indiquer l'apport du
texte par rapport \`a la litt\'erature existante.
}
au \S3
de commentaires bibliographiques.
Le \S2  d\'eduit le lemme d'\'elimination
d'une paire de points critiques en bonne position (le {\bf Corollaire C \/}en toute fin du \S1)
des {\bf Corollaires A'\/}, {\bf  B\/}, {\bf Lemmes C\/}, {\bf D\/}
et leurs compl\'ements.

Ces {\bf Corollaires A'\/}, {\bf A\/}  et {\bf B\/} suivent
des Lemmes de m\^eme lettre \'etablis, ainsi que leurs compl\'ements,
dans les quatre appendices A, B, C, D. correspondants.

Le trop laconique r\'esum\'e ci-dessus est d\'evelopp\'e en quatri\`eme de couverture.

Que Riccardo Benedetti et Sylvain Courte soient remerci\'es
pour leurs lectures critique du manuscript.

\npage 

\fancyhead[LO]{A. M. \quad C'est un petit val qui mousse de fonctions}
  \fancyhead[RO]{\thepage}
  \fancyhead[RE]{\S 1 Le fil des \'enonc\'es}
  \fancyhead[LE]{\thepage}
\centerline{\bf \S1 Introduction : d\'efinitions et fil des \'enonc\'es\/}
\vskip3mm 
\nDefns {\bf 0\/}
Soit une vari\'et\'e%
\footnote{\small
suppos\'ee s\'epar\'ee mais, avant le Lemme B p. 4,
non n\'ecessairement paracompacte.
}
lisse
$V$
munie d'une fonction r\'eelle lisse
$f :V\rightarrow{\Bbb R\/}$%
.

${\romannumeral 1})$
Un
{\it val
${\underline{\cal V\/}}_W$
pour\/}
$f$
est un microvoisinage tubulaire m\'etrique%
\footnote{\small
{\sl c. a d.\/}
un {\it voisinage tubulaire  ouvert\/}
${\cal V\/}_W$  
dans
$V$
de la sous-vari\'et\'e
$W$
de
$V$
muni d'un fibr\'e eucli\-dien
${\tilde \pi}\!: E\!\rightarrow\!W,\ {\tilde q}\!: E\!\rightarrow\!{\Bbb R\/}_{+}$
et d'une fonction de {\it taille\/}
$r_{f, W}\!=\!r\!:W\!\rightarrow]\!0, +\infty[$
positive lisse tels que, si
$E_{r}\!=\!\{x\!\in\!E\, ;\, {\tilde q}(x)\!\leq\!r({\tilde \pi}(x))\}$%
, il y a un diff\'eomorphisme 
$\tilde{\vartheta}\!:({\build{E_{2r}}_{}^{\circ}}, W)\!\rightarrow\!({\cal V\/}_W, W)$
avec
${\tilde q}\!=\!q\!\circ\!\tilde{\vartheta}$
et
${\tilde \pi}\!=\!\pi\!\circ\!\tilde{\vartheta}$%
, de restriction
$\vartheta\!=\!\tilde{\vartheta}_{|E_r}\!:E_{r}\!\rightarrow\!T$
dif\'eomorphisme sur le {\it tube ferm\'e\/}
$T$%
.}
$({\cal V\/}_W, T, \pi, q)$
%{\small
$$W\!\subset\!T\subset\!{\cal V\/}_W\!\subset\!V,\
\pi : {\cal V\/}_W\rightarrow W\
\hbox{\small  microfibr\'e euclidien
en boules de m\'etrique \/} q : {\cal V\/}_W\rightarrow{\Bbb R}$$
%}
d'une sous vari\'et\'e localement ferm\'ee
$W\subset V$
au voisinage de qui, pour une famille de diff\'eomorphismes
$\eta : W\rightarrow Diff_0({\Bbb R\/})$
 de
${\Bbb R\/}$
fixant des voisinages de l'origine,
l'application
$f$
est
$\eta$%
-translat\'ee de la
m\'etrique~:
$$f_{|T}-f_{|T}\circ\pi\!=\!(\eta\circ\pi)\cdot q:\quad \hbox{\rm pour tout \/}
x\in T, f(x)-f(\pi(x))=\eta(\pi(x))(q(x))$$
ainsi 
$f_{|\partial T}$
se factorise
$f_{|\partial T}=h_{f, W}\circ \pi$
par une fonction
$h_{f, W} : W\rightarrow{\Bbb R\/}$%
.
%$$\BoxedEPSF{val.eps }$$

On dira alors que 
$f$
est {\it fonction val\/}, de
${\romannumeral 1}0)\ $
{\it val \/}\
${\underline{\cal V\/}}_W$
et~:

 $${\romannumeral 1} 1)\,
\hbox{\it fond\/}\
W,\quad
{\romannumeral 1} 2)\,
\hbox{\it taille\/}\ 
r_W\!=\!r_{f, W},\ \hbox{\rm et \/}\quad
{\romannumeral 1} 3)\
\hbox{\it cr\^ete \/}
h_W\!=\!h_{f, W}$$

${\romannumeral 1} 4)$ 
Le val
${\underline{\cal V\/}}_W$
est {\it constant\/} pr\`es d'un compact
$K\subset W$
si son fibr\'e a un choix de microparam\'etrisations fibr\'ees%
\footnote{\small
coordinate bundle 2.3 de {\bf [St]\/} pour
$G={\bf O\/}({\mathop{\rm codim}\nolimits}_V(W))$
le groupe orthogonal.
}
$\bigl(\phi_j: W_j\times{\Bbb R\/}^r\rightarrow{\tilde{\pi}}^{-1}(W_j)\bigr)_{j\in J}$
dont le cocycle dans le groupe orthogonal
$g_{i, j}: W_i\cap W_j\rightarrow{\bf O\/}$
est constant pr\`es des
$K\cap W_i\cap W_j,\, i, j\in J$%
.

\Rmc {\bf Notations\/}|
Soit
$h, k : X\rightarrow{\Bbb R\/}, Y\!\subset\!X$
on note
$$h\!\leq_Y\!k\
\hbox{\rm (resp. \/}
h\!<_Y\!k
\hbox{\rm ) si pour tout \/}
y\in Y,\ h(y)\!\leq\!k(y)\
\hbox{\rm (resp. \/}
h(y)\!<\!k(y)
\hbox{\rm )\/}. $$

Les ensembles de fonctions lisses sur
$V$
et
%${\underline{\cal V\/}}_W$%
%-%
fonctions val de fond
$W$
sont not\'es~:
$${\goth A\/}(V)\!=\!{\goth A\/}\!=\!\{f : V \rightarrow {\Bbb R}\, ;\, f\ \hbox{\\small lisse\/}\}$$
$${\goth V\/}(V, W)\!=\!{\goth V\/}_W\!=\!\{f\!\in\!{\goth A\/}\, ;\, f\ \hbox{\\small est fonction val de fond \/} W\}$$

Sans  d\'efinir de structure de vari\'et\'e sur
${\goth A\/}(V), {\goth V\/}(V, W),\ldots$%
, une application
$\phi : {\cal D\/}\rightarrow{\goth V\/}(V, W)$
d\'efinie sur une partie
${\cal D\/}$
de
${\goth A\/}(V), {\goth V\/}(V, W)\!\times{\goth A\/}\!(W),\ldots$
est dite
{\it lisse\/}
si pour toute 
vari\'et\'e lisse
$X$
et
$\gamma : X\rightarrow{\cal D\/}$
d'\'evaluation lisse%
\footnote{\small
l'application
$(x, \underline{v})\!\in\!X\!\times\!V, X\!\times\!V\times\!W,\ldots\quad\mapsto\quad
{\overline \gamma\/}(x, \underline{v})\!=\!\gamma(x)(\underline{v})\in{\Bbb R}$
est lisse.
}
, alors l'application
$\overline{\phi\!\circ\!\gamma} : X\times V\rightarrow{\Bbb R\/},\
(x, v)\mapsto \overline{\phi\!\circ\!\gamma}(x, v)\!=\!\phi\circ\gamma(x)(v)$
est lisse.
\finc

\nDefns {\bf 1\/}
Soit
$a\!<\!b\!\in\!{\Bbb R\/}$%
, un {\it chemin de fonctions val de fond\/} 
$W$
[resp. {\it  fonctions\/}]
param\'etr\'e par l'intervalle
$[a, b]$
est 
$\gamma : {\Bbb R\/}\rightarrow{\goth V\/}(V, W)$%
[resp.
${\goth A\/}(V)$%
],
$s\mapsto \gamma_s$
lisse telle que pour tout
$s\!\leq\!a, \gamma_s\!=\!\gamma_a$
et pour tout
$s'\!\geq\!b, \gamma  _{s'}\!=\!\gamma_b$%
, les vals [resp. fonctions]
$\gamma_{-\infty}\!=\!\gamma_a$
et
$\gamma_{\infty}\!=\!\gamma_b$
sont dits {\it origine\/} et {\it extr\'emit\'e\/} du chemin
$\gamma$%
.
%\vfill\eject
\Thc Lemme A|
Soit
$(f, e)\!\in\!{\goth V\/}(V, W)\!\times\!{\goth A\/}(W)$%
, des compacts
$K, L\!\subset\!W$
avec~:
$K\!\subset\, \build{L}_{}^{\circ}\,\subset\!L,\
e\!\leq\!\frac{h_{ W}+2f_{|W}}{3},\ e-f_{|W}$
localement constante pr\`es de
$K$
et
$e_{|W\setminus\!\build{L}_{}^{\circ}}\!=\!f_{|W\setminus\!\build{L}_{}^{\circ}}$%
.%\hfill\break

Alors il y a un chemin
$f_{e, s}\!\in\!{\goth V\/}(V, W)$
lisse en
$(f, e, s)$
de
$f_{e, 0}=f$
\`a
$f_{e, 1}=f_{e}$
o\`u
${f_{e}}_{|W}\!=\!e$%
, pour tout
$s\!\in\!{\Bbb R\/}, h_{f_{e, s}, W}\!=\!h_{f, W},\
{f_{e, s}}_{|V\setminus(\build{T}_{}^{\circ}\cap\,{\pi_{|T}}^{-1}\!(\build{L}_{}^{\circ}))}\!=\!%
{f}_{|V\setminus(\build{T}_{}^{\circ}\cap\,{\pi_{|T}}^{-1}\!(\build{L}_{}^{\circ}))}$%
, pr\`es de 
$K, f_{e, s}-f$
 est localement constante et, si
 $e\!\leq\!f_{|W}$%
 , pour tout
$s\!\in\!{\Bbb R\/}, f_{e, s}\!\leq\!f$%
.
\finc

\nDefns {\bf 2\/}
Soit
$F\!\subset\!V, u'<u\!\in\!{\Bbb R\/}$
et une application
$f: V\rightarrow {\Bbb R\/}$%
.

Le {\it tronqu\'e de \/}
$F$
au dessus de
$u$
est
${_{f, u}}F\!=\!{_{u}}F\!=\!f^{-1}([u, +\infty[)\cap F$%
, l'application
$f$
est dite
$F$%
{\it -propre au dessus de\/}
$u'$
si %la restriction
$f_|\!: {_{u'}}F\!\rightarrow ]u', +\infty[$
est continue et propre%
\footnote{\small
c'est automatique si
$F$
est ferm\'e et l'application
$f$
est continue et propre.
}.%\hfill\break

Soit
$\kappa\!\in\!{\Bbb R\/}$%
, une {\it nappe (de
$C$
) pour $f$ au\/} niveau
$\kappa$
est le fond
$W\!\subset\!V$
d'un val de
$f$
avec
$\emptyset\!\ne\!{_{\kappa}}W\!=\!f^{-1}(\kappa)\cap W\!\build{=:}_{}^{\rm Def\/}\!C$
et 
$f_{|W\!\setminus\!C}$
est non singuli\`ere.

\Thc Corollaire  A|
Soit
$W$
une nappe
pour
$f$
au niveau
$\kappa$
et
$u'\!<\!u\!<\!\kappa$
tels que~:
$$f\ \hbox{\rm est \/}\ W\hbox{\rm -propre au dessus de \/} u'%
\leqno{\bigl({\cal P\/}{_{u'}}W\bigr)}$$
Alors si
${\cal U\/}$
est voisinage
de
${_{u}}W$
il y a un chemin 
\`a support dans
${\cal U\/}$
de fonctions
$f_{\sigma}$
pour lesquelles
$W$
reste nappe de
$C$
et un voisinage
${\cal C\/}$
de
$C$
avec si
$\sigma\!\leq\!\sigma'$%
~:
$${f_{\infty}}\!\leq\!f_{\sigma'}\!\leq\!f_{\sigma}\!\leq\!f_0\!=\!f,\quad
\hbox{\rm et\/}\quad f_{|{\cal C\/}}-\kappa+u\!=\!{f_{\infty}}_{|{\cal C\/}}$$
\finc
\vskip-2mm
\Rmc {\bf Notations\/}|
Si
$W$
est sous-vari\'et\'e d'une vari\'et\'e
$V$
et
$v\!\in\!V,\ w\!\in\!W$%
, on note~:%\hfill\break

${\goth A\/}_{v}$
l'anneau des  germes%
\footnote{\small
deux fonctions ont m\^eme germe en
$v$
si elles co\"{\i}ncident sur un voisinage de
$v$%
.
}
en
$v$
de l'anneau
${\goth A\/}\!=\!{\goth A\/}(V)$
des fonctions lisses sur
$V$%
.

${\goth W\/}\!=\!\{a\!\in\!{\goth A\/}\,;\,a_{|W}\!=\!0\}\!\subset\!{\goth M\/}\!=\!{\goth M\/}(w)\!=\!\{b\!\in\!{\goth A\/}\,;\, b(w)\!=\!0\}\!\subset\!{\goth A\/}$
les id\'eaux nuls sur
$W$%
et en
$w$
respectivement et
${\goth W\/}_w\!\subset\!{\goth M\/}_w\!\subset\!{\goth A\/}_w$
les id\'eaux de germes correspondants.

L'anneau quotient
${\goth B\/}\!=\!{\goth A\/}/{\goth W\/}$
s'identifie \`a l'anneau des fonctions lisses sur
$W$%
.

On note
$\rho : {\goth A\/}\rightarrow{\goth B\/},\ \rho_{w}: {\goth A\/}_{w}\rightarrow{\goth B\/}_{w}$
les morphismes quotient de restriction.

Le
${\goth A\/}$%
-module
${\cal X\/}(V)$
 des champs de vecteurs%
\footnote{\small
${\cal X\/}(V)\!=\!\{X : {\goth A\/}\rightarrow{\goth A\/}\,;\,\hbox{\rm si \/}
f, g\!\in\!{\goth A\/},\, X(f+g)\!=\!X(f)+X(g),\,
X(f\cdot g)\!=\!X(f)\cdot g+f\cdot X(g)\}$
}
sur
$V$
 contient le sous-module
${\cal Y\/}\!=\!{\cal Y\/}_{W}(V)\!=\!%
\{Y\!\in\!{\cal X\/}\,;\, Y({\goth W\/})\!\subset\!{\goth W\/}\}$
des  {\it champs tangents \`a \/}
$W$%
.

Pour
$v\!\in\!V, w\!\in\!W,\, Z\!\in\!{\cal X\/}, Y\!\in\!{\cal Y\/}$%
, on note
$Z_{v}\!:{\goth A\/}_{v}\!\rightarrow\!{\Bbb R\/},%
Y_{w}\!:{\goth B\/}_{w}\!\rightarrow\!{\Bbb R\/}$
les {\it vecteurs tangents \`a\/}
$V$
et
$W$
en
$v$
et
$w$%
, d\'efinis par
$Z_{v}(a)\!=\!Z(a)(v)$
et
$Y_{w}(b)\!=\!Z(b)(w)$%
, \'el\'ements
$Z_{v}\!\in\!{\mathop{\rm T\/}\nolimits}_{v}V,
Y_{w}\!\in\!{\mathop{\rm T\/}\nolimits}_{w}W$
des espaces tangents%
\footnote{\small
o\`u
$T_{w}V\!:=\!\{\xi: {\goth A\/}_{v}\!\rightarrow\!{\bf R\/};
{\bf R\/}\hbox{\rm -lin\'eaires et si
}
f_{w}, g_{w}\!\in\!{\goth A\/}_{v},
\xi(f_{w}\cdot g_{w})\!=\!\xi(f_{w})\cdot g_{w} + f_{w}\cdot \xi(g_{w})\},\ldots$%
. De m\^eme si
$\phi: U\rightarrow U'$
est une application diff\'erentiable entre deux vari\'et\'es lisses
$U, U'$%
, son application tangente en
$u\in U$
est
$T\phi_{u}: T_{u}U\rightarrow T_{\phi(u)}U',\ T\phi_{u}(\xi)(b)=\xi(b\circ\phi_{u})$%
.
} \`a
$V$
et
$W$
en
$v$
et
$w$
corres\-pondants qui s'identifient%
\footnote{\small
car
${\goth W\/}{\cal X\/}\!\subset\!{\cal Y\/}$%
, m\^{e}me si
${\goth w\/}\!\in{\goth W\/}, X\!\in\!{\cal X\/},
{\goth w\/}X({\goth A\/})\!\subset\!{\goth W\/}$%
, car
${\goth W\/}$
\'etant id\'eal \`a droite,
${\goth w\/}X(\alpha)\!\in\!{\goth W\/}$
et notant, si
$v\!\in\!V, w\!\in\!W, {\cal X\/}_v, {\cal Y\/}_w$
 les
${\goth A\/}_v$
-modules de germes de champs de vecteurs associ\'es.
.
}aux quotients%
~:
$${\mathop{\rm T}\nolimits}_{w}V\!=\!{{\cal X\/}_w/}_{{\goth M\/}_w{\cal X\/}_w},\
{\mathop{\rm T}\nolimits}_{w}W\!=\!{{\cal Y\/}_w/}_{({\goth M\/}_w{\cal Y\/}_w+{\goth W\/}_w{\cal X\/}_w)}$$

L'{\it id\'eal Jacobien de
$f$
\/}
est~:
${\goth J\/}(f)\!=\!{\cal X\/}(f)\!=\!\{X(f)\,;\,
X\!\in\!{\cal X\/}\}\!\subset\!{\goth A\/}$%
, celui de la restriction
$f_{|W}$
de
$f$
\`a
$W$
est donc
${\goth J\/}(f_{|W})\!=\!\rho({\cal Y\/}(f))\!=\!
\{\rho(Y(f); Y\!\in\!{\cal Y\/}\}\!\subset\!{\goth B\/}$%
.

\finc

\nDefns {\bf 3\/}
Soit
$W\!\subset\!V$
une sous-vari\'et\'e et une fonction 
$f: V\!\rightarrow\!{\Bbb R\/}$
lisse.

${\romannumeral 1}0)$
Un champ de vecteur
$Z$
d\'efini pr\`es de
$W$
est
$f$%
{\it -plat\/} si
$Z(f)\!\in\!{\goth W\/}$%
.

${\romannumeral 1}1)$
Un tube
$\pi: T\rightarrow W$
autour de
$W$
dans
$V$
est
$f$%
-{\it plat\/}
si
$Ker(T\pi)\subset ker(Tf)$
{\it c.ad.\/} tout champ de vecteurs
$Z$
sur
$T$
port\'e%
\footnote{\small
{\it c. a d.\/} pour tout
$t\!\in\!{\mathop{\cal N\/}_{W}\nolimits},
{\mathop{\rm T\/}\nolimits}_{t}(\pi)(Z_{t})\!=\!0%
\!\in\!{\mathop{\rm T\/}\nolimits}_{\pi(t)}W$%
.
}
par les fibres de
$\pi$
est
$f$%
-plat.

${\romannumeral 1}2)$
La sous-vari\'et\'e
$W$
est
$f$%
-{\it directe dans \/}
$V$
si~:
$$\hbox{\rm pour tout \/} w\!\in\!W,\
\rho_{w}({\goth J\/}_{w}(f))\subset{\goth J\/}_{w}(f_{|W})\leqno({\cal D\/})$$

${\romannumeral 2}1)$
Les {\it ensemble  critique\/}
et {\it ensemble critique sur\/}
$W$
de
$f$
sont\footnote{\small
${\cal X\/}_c(f)(c)\!=\!0$
 signifiant
 pour tout 
$X\!\in\!{\cal X\/}_c,\ X(f)(c)\!=\!0$% 
.}~:
$${\mathop{\Sigma}\nolimits}(f)\!=\!\{c\in V\, ;\, {\cal X\/}_c(f)(c)\!=\!0\}%
\quad
{\rm et\/}\quad
{\mathop{\Sigma}\nolimits}_W(f)\!=\!\{d\in W\, ;\, {\cal Y\/}_d(f)(d)\!=\!0\}\!=\!%
{\mathop{\Sigma}\nolimits}(f_{|W})$$

%, tout point critique de
%$f_{|W}$
%est point critique de
%$f$%

${\romannumeral 2}2)$
Un point
$w\!\in\!{\mathop{\Sigma}\nolimits}_W(f)$
critique de la sous-vari\'et\'e
$W$
est 
${\goth W\/}$%
{\it -d\'efini positif\/} si~:
%pour la Hessienne%
%
%
$${\goth W\/}_{w}\subset{\goth J\/}_{w}(f)\subset{\goth M\/}_{w}
\leqno(D_{\goth W\/})$$
$$\hbox{\rm pour tout \/} X\!\in\!{\cal X\/}_{w}\
f\hbox{\rm -plat\/},
\ X(X(f))(w)\!\geq\!0 \leqno(P)$$
{\small
[Si  
$f_|W$
est de Morse%
\footnote{\small
voir les {\bf D\'efinitions 4\/}
${\romannumeral 3)}, {\romannumeral 1)}$
ci dessous et la fin de l'appendice B%
. 
}
 traduisez~:
 si
 $y\!\!\in\!\!\Sigma_W(f), y\!\!\in\!\!\Sigma(f)$
 non d\'eg\'en\'er\'e avec
 ${\mathop{\rm Ind}\nolimits}_y(f)\!\!=\!\!{\mathop{\rm Ind}\nolimits}_y(f_{|W})$%
.]
}

${\romannumeral 3}1)$
$f$
est {\it stabilis\'ee en\/}
$W\!$%
, ou
$W$%
{\it -stabilis\'ee\/}, si 
elle est fonction val de fond
$W$%
.

${\romannumeral 3}2)$
$f$
est {\it infinit\'esimalement stabilis\'ee\/}, ou
${\goth W\/}$%
{\it -stabilis\'ee,\/} en
$W$
si
$W$
est\break
$f$%
-directe dans
$V$
et tout
$w\!\in\!{\mathop{\Sigma}\nolimits}_W(f)$
point critique sur
$W$
est
${\goth W\/}$%
-d\'efini positif.

\Thc Lemme B |
Soit
$W\!\subset\!V$
une sous-vari\'et\'e paracompacte.

$(1)$
$W$
 a un tube
$f$%
-plat si et seulement elle est
$f$%
-directe dans
$V$%
.

$(2)$
$f$
est stabilis\'ee en
$W$
si et seulement si elle l'est
%${\mathop{\Sigma}\nolimits}_W$%
infinit\'esimalement.% stabilis\'ee.

\finc
Dans la suite toutes les vari\'et\'es sont suppos\'ees paracompactes.

{\small%\parc
Appliquant le {\bf Lemme B\/} \`a
$f$
et
$W\!\subset\!V$
contenant un point
$c\!\in\!{\mathop{\Sigma}\nolimits}(f)\cap W$
critique de
$f$
avec
$N\!:=\!{\mathop{\rm T}\nolimits}_cW\!\subset\!{\mathop{\rm T}\nolimits}_cV$
d\'efini n\'egatif pour
$q_c$
de dimension maximale,
puis \`a
$-f_{|W}$
et
$\{c\}\!\subset\!W$%
~:
}%\finc
\Thp Corollaire B
{\small\rm (Lemme de Morse de forme canonique
des points critiques non d\'eg\'en\'er\'es)}|
Soit
$c\in{\mathop{\Sigma}\nolimits}(f)$
point critique non d\'eg\'en\'er\'e%
\footnote{\small
{\it c.ad.\/} la d\'eriv\'ee seconde
$b_{c}$
de
$c$
induit un isomorphisme de
${\mathop{\rm T\/}\nolimits}_{c}V$
sur son dual
${\mathop{\rm T\/}\nolimits}_{c}^{\ast}V$%
.
}
de
$f\!: V\!\rightarrow\!{\Bbb R\/}$
et 
$N\!\subset\!E\!=\!{\mathop{\rm T}\nolimits}_cV$%
\hfill\break
(resp.
$P\!\subset\!E$%
)
sous espace d\'efini n\'egatif (resp. positif) maximal
pour la Hessienne%
\footnote{\small
rappelons que si
$w\!\in\!{\mathop{\Sigma}\nolimits}(f)$
la {\it forme quadratique Hessienne de 
$f$
en \/}
$w$
(terme quadratique de la formule de Taylor en
$w$%
)
est
${q_H}_w\!=\!q_w : {\cal X\/}_w\rightarrow{\bf R\/},\
{q_H}_w(X)\!=\!\frac{1}{2}X(X(f))(w)$
se factorise par l'espace tangent
${\mathop{\rm T}\nolimits}_wV\!=\!%
{{\cal X\/}_w/}_{{\goth M\/}_w{\cal X\/}_w}$
et a pour forme bilin\'eaire sym\'etrique associ\'ee la
{\it d\'eriv\'ee seconde \/}
$b_w\!=\!b : {\cal X\/}_w\times{\cal X\/}_w\rightarrow{\bf R\/},\
b(X, X')\!=\!q(X+X')-q(X)-q(X')\!=\!X(X'(f))(w)$%
.
}\break
$q\!=\!q_c : {\mathop{\rm T}\nolimits}_cV\rightarrow{\Bbb R\/}$
de
$f$
en
$c$
et
$P=N^{\perp_b}$
(resp.
$N=P^{\perp_b}$%
)
son 
$b$%
-orthogonal.

Alors
$0$
a des voisinages \'etoil\'es
$U_N, U_P,\, U\!=\!U_N+U_P$
dans
$N, P,\ E\!=\!N\!\oplus\!P$
et un diff\'eo\-%break%
morphisme
de param\'etrisation 
$\psi : (U, 0)\!\rightarrow\!({\cal U\/}, c)\subset(V, c)$
avec~:
$${\mathop{\rm T}\nolimits}_0\psi\!=\!{\mathop{\rm Id}\nolimits}_{{\mathop{\rm T}\nolimits}_cV}\quad
\hbox{\rm et\/}\quad
f\circ\psi-f(c)=q_{|U}=(q_{|N}\oplus q_{|P})_{|U}=(-(-q_{|N})\oplus q_{|P})_{|U}$$
\finp
\vskip-5mm
\rmc {\bf Remarque\/}|
Si
$N'\!\subset\!E$
est un autre sous-espace d\'efini  n\'egatif maximal alors
$N\!\cap\!P=0$
et
$N'=\{(x, g_{N, N'}(x))\!\in\!N\!\oplus\!P\!=\!E\, ;\, x\!\in\!N \}$
o\`u
$g_{N, N'}\!\in\!{\mathop{\rm Hom }\nolimits}(N, P)$%
. L'ensemble des choix de
$N$
s'identifiant \`a
$\{g\!\in\!{\mathop{\rm Hom }\nolimits}(N, P)\, ;\,  q_{|P}\circ g<_N-q_{|N}  \}$%
, convexe dans
${\mathop{\rm Hom }\nolimits}(N, P)$%
, est contractile.
\finc
\nDefns {\bf 4\/}
${\romannumeral 1})$
L'{\it indice\/}
et le {\it coindice\/} en le point 
$c\!\in\!{\mathop{\Sigma}\nolimits}(f)$
critique de
$f$%
, not\'es respectivement
${\mathop{\rm Ind\/}\nolimits}_{c}(f)$
et
${\mathop{\rm Ind\/}\nolimits}_{c}(-f)$%
, sont les dimensions de sous-espace respectivement
$q_c$
n\'egatif ou nul et
$q_c$
positif ou nul maximaux de
${\mathop{\rm T}\nolimits}_cV$%
.

${\romannumeral 2})$
Un voisinage d'un point critique non d\'eg\'en\'er\'e
est {\it de Morse\/} si il a une param\'etrisation, dite
{\it de Morse\/},
donn\'ee par le {\petcap Corollaire B\/}.

${\romannumeral 3})$ La fonction 
$f$
est {\it de Morse\/} si tous ses points critiques sont non d\'eg\'en\'er\'es.

${\romannumeral 4})$
Une {\it m\'etrique adapt\'ee \`a\/}
$f$
est une m\'etrique qui dans
des param\'e\-trisations de Morse de ses points critiques est
$\mu\!=\!\mu^N\oplus\mu^P : {\mathop{\rm T}\nolimits}_cV\!=\!N\oplus P\!\rightarrow\!{\Bbb R\/}$
o\`u
$\mu^N, \mu^P$
sont les formes quadratiques d\'efinies positives~:
$\mu^N=-q_{|N} : N\!\rightarrow\!{\Bbb R\/}, \mu^P=q_{|P} : P\!\rightarrow\!{\Bbb R\/}$%
.

${\romannumeral 5})$
Un champ de vecteurs 
$Z$%
, sur
$X$
est {\it pseudo-gradient faible\/} (abr\'eg\'e en {\it  pgf\/})
d'une fonction lisse
$f$
si en tout point r\'egulier
$x\!\in V\!\setminus\!{\mathop{\Sigma}\nolimits}(f)$
on a
$Z(f)(x)\!>\!0$%
. %\hfill\break

Si de plus
$f$
est de Morse et dans des para\-m\'e\-tri\-sations de Morse de ses points critiques,
le pgf
$Z\!=\!{\mathop{\rm grad}\nolimits}_{\mu}(f)$
est  gradient de
$f$
pour une m\'etrique adapt\'ee,\break
ce champ de vecteur
$Z$ est dit {\it pseudo-gradient\/} (abr\'eg\'e en {\it  pg\/}) de
$f$%
.

${\romannumeral 6})$
Une {\it donn\'ee\/}
${\cal F\/}\!=\!(f,  Z)$
[resp. {\it donn\'ee de Morse\/}
${\cal F\/}\!=\!(f, \Psi, Z)\build{=}_{}^{{\rm abr.\/}}(f,  Z)$%
]
est une fonction [resp. fonction de Morse]
$f$
munie d'un pgf %pseudo-gradient faible
[resp. pg%pseudo-gradient
]
$Z$
[resp. et d'une famille
$\Psi\!=\!(\psi_c : (U_c,0)\!\rightarrow\!({\cal U\/}_c, c)\bigr)_{c\in{\mathop{\Sigma}\nolimits}(f)}$
de para\-m\'e\-trisations de Morse dans lesquelles
$Z\!=\!{\mathop{\rm grad}\nolimits}_{\mu}(f)$
est, pour une m\'etrique adapt\'ee
$\mu$%
, gradient de
$f$%
.]

${\romannumeral 7})\
{\cal F\/}\!=\!(f, Z)$
est
$Z'$%
-\nobreak{\it modifi\'ee\/} de %'une donn\'ee
%(resp. donn\'ee de Morse)
${\cal F'\/}\!=\!(f', Z')$
, not\'e
${\cal F\/}\!\preceq\!{\cal F'\/}$%  
, si
$Z\!=\!Z'$%
\footnote{\small 
  et donc
  $f$
  et
  $f'$
  ont m\^eme ensemble ciritique
  $\Sigma(f)\!=\!\Sigma(f')$%
  .
}%
.
%et~:
%$$\hbox{\rm pour tout \/} c\in\Sigma(g'),\  U_c\!\subset\!U'_c,\ V_c\subset V'_c,\
%\psi'_c : U'_c\rightarrow V'_c\ \hbox{\rm \'etend \/} \psi_c={\psi'_c}_{|U_c}$$

${\romannumeral 8})$
Soit
$W\!\subset\!V$
 sous-vari\'et\'e directe sur qui la restriction
$f_{|W}$
est de Morse.
Une {\it donn\'ee de Morse adapt\'ee \`a \/}
$W$
est une donn\'ee de Morse%
\footnote{\small
l'existence d'une donn\'ee de Morse
${\cal F\/}\!=\!(f, \Psi, Z)$
sous-entend que la fonction
$f$
est de Morse.\hfill\break
La sous-vari\'et\'e \'etant directe, les points singuliers
de
$f_{|W}$
sont donc des points de Morse de
$f$%
.
}
${\cal F\/}\!=\!(f, \Psi, Z)$
dont
%${\romannumeral 8}1)$
le champ
$Z$
est tangent \`a
$W$
et
%${\romannumeral 8}2)$
si
$c\!\in\!{\mathop{\Sigma}\nolimits}(f_{|W})$%
, la param\'etrisation de Morse
$\psi_c$
de
${\cal U\/}_c$
induit une param\'etrisation de Morse
${\psi_c}_{|U_c\cap{\mathop{\rm T}\nolimits}_c(W)}$
du voisinage
 ${\cal U\/}_{W, c}\!=\!{\cal U\/}_c\cap W$%
.

${\romannumeral 9})$
Un val
${\underline{\cal V\/}}_W$
pour une fonction 
$f$
ayant une donn\'ee de Morse adapt\'ee \`a
$W$
est {\it de Morse\/}
si son tube a en tout
$c\!\in\!{\mathop{\Sigma}\nolimits}(f_{|W})$
pour (micro)param\'etrisation fibr\'ee, une donn\'ee de Morse adapt\'ee \`a
$W$
d'un voisinage
${\cal U\/}_c$
de
$c$%
.

%${\romannumeral 7} 4 )$
%Si
%$u\!<\!g(c)$
%les {\it tronqu\'ees en $u$ \/}  des nappes sont~:
%$$(\overline{N}_{u}, N_{u})=\bigl(\overline{c}_{_{\searrow\!\cdots\!\searrow}}\!(Z)_u, c_{_{\searrow}}(Z)_%u\bigr)=
%\bigl(\overline{c}_{_{\searrow\!\cdots\!\searrow}}\!\!(Z), c_{_{\searrow}}\!(Z)\bigr)\cap g^{-1}([u, +\inf%ty[)$$
%$$(\overline{P}^{v}, P^{v})=\bigl({^{^\nwarrow\!\cdots\!\nwarrow}\overline{c}}(Z)^v, {^{^\nwarrow}c}(Z)^v\%bigr)=
%\bigl({^{^\nwarrow\!\cdots\!\nwarrow}\overline{c}}(Z), {^{^\nwarrow}c}(Z)\bigr)\cap g^{-1}(]-\infty, v])$$

\Thc Compl\'ement B|
Si 
$f$
a une donn\'ee de Morse adapt\'ee \`a
$W$
et  est in\-fi\-ni\-t\'esimalement stabilis\'ee,
elle est stabilis\'ee en
$W$
et munie d'un val de Morse
${\underline{\cal V\/}}_W$%
.
\finc
\vskip-2mm
\Thc Compl\'ement A|
Si de plus
$e\!\leq\!f_{|W}, {\mathop{\Sigma}\nolimits}_W(f)\!\subset\!K$
et
$f$
a une donn\'ee 
${\cal F\/}\!=\!(f, Z)$
avec
$Z_{|W}$
pseudo-gradient faible (resp. pseudo-gradiant)
de
$e$
alors~: %\hfill\break

Quitte \`a r\'eduire la taille du val
${\underline{\cal V\/}}_W$%
, le champ
$Z$
est pgf (resp. pg) de
$f_{e}$%
.
\finc

\nDefns {\bf 5\/}
Notant
$\Phi_t$
le flot%
\footnote{\small
$\Phi : V\times{\bf R\/}\supset{\cal Z\/}\rightarrow V, (x, t)\mapsto z_x(t)\
\hbox{\rm solution maximale de \/}
z(0)\!=\!x\
\hbox{\rm et \/}
z'(t)\!=\!Z(z(t))$%
.
}
de
$Z$%
, pour
$x, y\!\in\!V$
deux points de
$V$
 on dira~:

${\romannumeral 1})$
$x$
est {\it ($Z$-)avant\/}
$y$
(not\'e
$x\!\leq^Z\!y$%
) si il y a
$b\in[0, +\infty]$
avec
$\build{Lim}_{t\rightarrow b}^{}{\Phi_t(x)}=y$%
.

${\romannumeral 2})$
$x$
est {\it ($Z$-)sous\/}
$y$
(not\'e
$x\!\preceq^Z\!y$
) si
$x$
est, pour un 
$n\!\in\!{\Bbb N\/} $%
,
{\it $nZ$-sous\/}
$y$%
, d\'efini r\'ecursivement par : soit
$n\!=\!0$
et
$x\!\leq\!^{Z}y$%
, soit 
il y
$y_{1}\!\leq\!^{Z}y$
avec
$x\!\preceq_{n-1}^{Z}\!\build{Lim}_{t\rightarrow-\infty}^{}\Phi_{t}(y_{1})$%
,
Les relations
$\leq^{Z}, \preceq^{Z}$
sont des ordres sur
$V$%
.
Soit
$c\!\in\!{\mathop{\Sigma}\nolimits}(f)$
un point critique de
$f$%
.

${\romannumeral 3})$
La {\it nappe de\/}
$c$
par
$Z$
est
${c_{_\searrow}}\!(Z)=\,\{x\in V\,;\, x\leq^{Z} c\}$%
.

${\romannumeral 4})$
La {\it conappe de\/} 
$c$
{\it par\/}
$Z$
est 
${^{^\nwarrow}c}(Z)\!=\!c_{_\searrow}\!(-Z)$%
.

${\romannumeral 5})$
La {\it  nappe satur\'ee de\/}
$c$
{\it par\/}
$Z$
est le ferm\'e
$\overline{c}_{_{\searrow\!\cdots\!\searrow}}\!\!(Z)=\,\{x\in V\,;\, x\preceq^{Z}c\}$

${\romannumeral 6})$
Le {satur\'e de\/}
$c$
{\it par\/}
$Z$
est le ferm\'e
$\overline{c}(Z)\!=\!\overline{c}_{_{\searrow\!\cdots\!\searrow}}\!\!(Z)\cup\overline{c}_{_{\searrow\!\cdots\!\searrow}}\!\!(-Z)$
\rmc Remarque A|
$(1)$
Si
$(f, Z)$
est une donn\'ee de Morse la nappe
${c_{_\searrow}}\!(Z)\!=\!W$
d\'efinie en
${\romannumeral 3})$\hfill\break
est une sous-vari\'et\'e, nappe pour la D\'efinition 2.
Laissons la lectrice traduire le Corollaire A~:
\finc
\rmc {\bf Corollaire A'\/} {\rm (lemme d'abaissement de Morse)\/}|
{\sl Si
$(f, Z)$
est une donn\'ee de Morse, le chemin
$f_s$
du {\bf Corollaire A\/}
produit un chemin
$(f_s, Z)\!\preceq(f, Z)$
de donn\'ees de Morse
$Z$% 
-modifi\'ees.
}
\finc
\vskip-2mm
\rmc $(2)$|
Si, au lieu de la condition
$\bigl({\cal P\/}{_{u'}}c\bigr)\build{=}_{}^{\rm Def\/}\bigl({\cal P\/}{_{u'}}W\bigr)$%
du {\bf Corollaire A'\/}, on a%
\footnote{\small
une condition plus faible, car automatiquement satisfaite si
$f$
est propre.
}~:
$$f\ \hbox{\rm est \/}
\overline{c}_{_{\searrow\!\cdots\!\searrow}}\!(Z)%
\hbox{\rm -propre au dessus de \/} u' \leqno{\bigl(\overline{\cal P\/}{_{u'}}c\bigr)}$$
la nappe satur\'ee tronqu\'ee
${_{u}}\overline{c}_{_{\searrow\!\cdots\!\searrow}}\!(Z)$
n'a, hormis
$c$%
, qu'un nombre fini de points critiques. %\hfill\break

Une
$Z$%
-modification laissant inchang\'ees nappes et nappes satur\'ees, le {\bf Corollaire A'\/}
conduit%
\footnote{\small
par induction croissante  sur l'ordre
$\preceq_Z$
de ces points critiques.
}\
 au  m\^eme \'enonc\'e en rempla\c cant l'ouvert
${\cal U\/}\supset {_{u}}c_{_{\searrow}}\!(Z)$
par
${\cal V\/}\supset {_{u}}\overline{c}_{_{\searrow\!\cdots\!\searrow}}\!\!(Z)$%
.
\finc

\rmc $(3)$|
En appliquant {\bf Corollaire A'\/}
et remarque A
$(2)$
\`a la fonction oppos\'ee
$-f$
on obtient pour
$f$
deux lemmes d'\'el\'evation d'\'enonc\'es laiss\'es aussi aux bons soins de la lectrice, puis au~:
\finc

\rmc {\bf coCorollaire A'\/} {\rm (lemme de d\'eplacement de Morse)\/}|
{\sl Soit une donn\'ee de Morse
$(f, Z)$
et
$g  : (\Sigma(f), \preceq^{Z})\rightarrow({\bf R\/}, \leq)$
strictement croissante,
\'egale \`a
$f_{|\Sigma}$
hors d'un ensemble fini.\hfill\break
Alors si pour tout
$c\!\in\!\Sigma(f)$
avec
$g(c)\!\ne\!f(c)$
la restriction
$f_{|\overline{c}(Z)}$
est propre et
${\cal U\/}$
est voisinage de
$$\bigcup_{c\in\Sigma, g(c)\ne f(c)}\overline{c}(Z)\cap f^{-1}([\min(f(c), g(c)), \max(f(c), g(c))])$$
il y a un chemin lisse
$(f_s, Z)\preceq(f, Z)$
\`a support dans
${\cal U\/}$
de donn\'ees de Morse
$Z$%
-modifi\'ees de
$f_{-\infty}\!=\!f$
\`a
$f_1\!=\!f_{\infty}$
tel que, pr\`es de tout
$c\!\in\!\Sigma, f_s\!=\!f+f_s(c)-f(c)$
et
${f_1}_{|\Sigma}\!=\!g$%
.\hfill\findem
}
\finc

\vskip5mm

\parc
Avec les notations
$c\in\Sigma(f), N, P\!\subset\!E, \psi :U\!=\!U_N\!+\!U_P\!\rightarrow\!{\cal U\/},\ q$
du Corollaire B,
$\mu\!=\!\mu^N\!\!\oplus\!\mu^P\!\!, Z$
de la {\bf D\'efinition 4\/},
le champ
$Z^{U}\!\!\build{=}_{}^{\rm Def\/}\!\psi^{-1}_{\ast}(Z_{|{\cal U}})$%
, de flot
$\Phi^{U}_t$%\
, restriction \`a
$U$
de %u flot lin\'eaire
$\Lambda_t\!=\!e^{-t{\mathop{\rm Id}\nolimits}_N+t{\mathop{\rm Id}\nolimits}_P}$%
, a pour orbites les intersections avec
$U$
de celles de
$\Lambda_t$%
.
Soit
$p\!\in\!P,\ \mu^{P}(p)\!=\!1$
et
$R\!=\!p^{\perp_{\!\mu^P}}$%
, notons
$\gamma\!=\!{\bf R\/}_+^{\ast}\!\cdot\!p\cap U,\
{\cal U\/}^{\epsilon}\!=\!\{u\in\!E\,;\, \mu(u)\!<\!\epsilon^2\},\ 
(N^{\epsilon}, P^{\epsilon}, R^{\epsilon})\!=\!(N, P, R)\cap {\cal U\/}^{\epsilon}$%
. Si
$
s\!>\!0$
avec
$s\cdot p\!\in\!\gamma$
on a~:
\finc

\Thp Lemme C {\rm\small (extension canonique d'une nappe transverse)}|
Soit 
$M\!\subset\!U$
une sous-vari\'et\'e invariante%
\footnote{\small
Une partie
$A\!\subset\!U$
est dite {\it invariante par\/}
$\Phi^U$
si pour tout
$u\!\in\!U$
et
$t\!\in\!{\Bbb R\/}$
telle que
$\Phi^U_t(u)$
est d\'efini alors
$a\!\in\!A$
si et seulement si
$\Phi^U_t(a)\!\in\!A$%
.
}
par le flot
$\Phi_t^{U}$%
, 
et telle que le rayon
$\gamma$
est ouvert dans
$M\cap P$
avec intersection  transverse.
Alors~:

$(1)$ Il y a
$s, \nu, \rho\!>\!0,\ N^{\nu}\!\!\subset\!U_N,\, s\cdot p+R^{\rho}\!\subset\!U_P$
et
$\theta\!:\!(N^{\nu}, 0)\!\rightarrow\!(R, 0)$
lisse avec~:
$$\bigl(M-s\cdot p\bigr)\cap\bigl(N^{\nu}+R^{\rho}\bigr)\!=\!\{n+\theta(n)\, ;\, n\in N^{\nu}\}$$

$(2)$
Il y a
$0\!<\!\delta$ %\!\leq\!\sqrt{s\nu}$
tel que l'application~:
$$\chi : \bigl(N^{\delta}, \{0\}\bigr)\times]-\delta, \delta[\rightarrow E,\ \
\chi(n, \tau)=n+\tau p+\frac{\tau}{s}\!\cdot\!\theta(\frac{\tau}{s}\!\cdot\!n)$$
est une paire de plongements d'image tangente au champ
$Z^{U}$
avec~:
$$\chi(0)\!=\!0,\quad \chi(\{0\}\times]0, \delta[)=\gamma\cap{\cal U\/}^{\delta}\quad
\hbox{\rm et\/}\quad
\chi(N^{\delta}\times]0, \delta[)\subset M\cap\{q\!>\!0\}\!:=M'$$

$(3)$
Le seul point critique, de
$f\circ\psi\circ\chi$
 est
 $0$%
 , il est non d\'eg\'en\'er\'e de coindice
 $1$%
 .
\finp
% \vskip7mm
%\begin{figure}
%\centerline{\BoxedEPSF{extension.eps scaled 250 } }
\TrimTop{-10pct}
\TrimBottom{-7pct}
\centerline{\BoxedEPSF{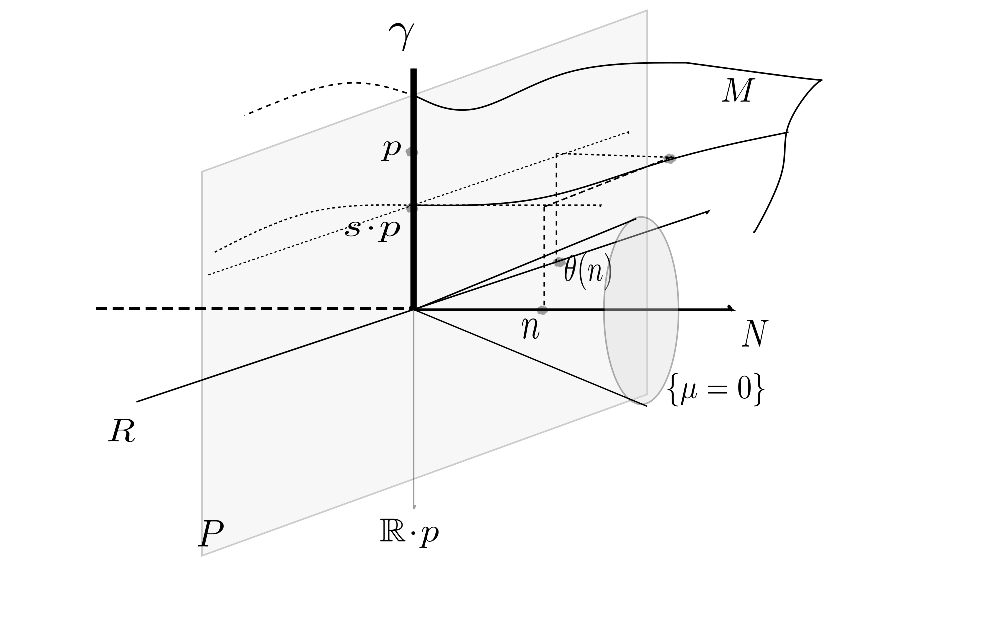 scaled 500} }
%\Lafigps{extension.eps}{500}{extension de nappes}
%\end{figure}
\vskip4mm
\centerline{Param\'etrisation d'une sous-vari\'et\'e invariante par le flot
$\Phi_{t}^{U}$}
\vskip5mm
%hfill\break
%{\small
Ainsi il y a une paire de vari\'et\'es
$(M_{\chi}, l_{\chi})$
telle que
$c\in l_{\chi}\!\subset\!M_{\chi}$
et contenant
$(M', \gamma)\!\subset\!(M_{\chi}, l_{\chi})$
comme sous-paire ouverte et
une immersion
$\xi : M_{\chi}\rightarrow U$
avec
$c$
unique point critique des restrictions
$f\circ\xi_{|M_\chi}, f\circ\chi_{|l_\chi}$%
, les  coindices \'etant
$1$%
.
%}

\Thc Compl\'ement C| 
Il y a un voisinage
$(M^{+}, \lambda^{+})$
de
$(\{0\}\!\cup\!M', \{0\}\!\cup\!\gamma)$
dans 
$(M_{\chi}, l_{\chi})$
sur lequel 
$\xi, \xi_{|l_\chi}$
induisent des plongements si et seulement si~:

$(4)$
Le rayon
$\gamma=P\cap M$
est toute l'intersection de
$P$
avec
$M$%
.
 
En ce cas %leurs images
$\psi\circ\xi(M^{+}), \psi\circ\xi(\lambda^{+})$
sont fonds de vaux de Morse pour
$f$
et
$-f_{|M^{+}}$%
.
\finc

\Thp Lemme D {\rm\small
  (Lemme du dromadaire)}| 
Soit
$I\!\subset\!{\Bbb R\/}$
un intervalle ouvert,
$k :I\!\rightarrow\!{\Bbb R\/},\ b\!<c\!<\!d\!<\!n\in\!I$
avec~:
$${k'}^{-1}(]-\infty, 0], \{0\})=\bigl([c, d], \{c, d\}\bigr)\
\hbox{\rm et \/}\ k(c)<k(n)\leqno(D)$$

Alors il y a un chemin
$k_t$%
, 
issu de
$k_0\!=\!k$
\`a
$k_{\infty}\!=\!k_2$
et
$t_0\in\, ]1, 2[$
avec~:

$(1)$
pour tout
$t\!\leq\!t',\ k_t\!\leq\!k_{t'}$
avec \'egalit\'e si
$x\not\in\, ]b, n[$%
.

$(2)$
Si
$0\!<\!t\!<\!t_0,\ k_t$
est de Morse,
et si
$t_0\!<\!t,\
{k_t'}^{-1}(]-\infty, 0])\!=\!\emptyset$%
, de plus~:

$(3)$
pr\`es de
$[c, d]$
et pour
$1\!\leq\!t,\ k_t$
est restriction d'un polyn\^ome cubique.

\finp
%\vskip-4mm

%%%
%\midinsert
\TrimTop{-10pct}
\TrimBottom{-5pct}
\TrimLeft{-5pct}
\TrimRight{-5pct}
\centerline{\BoxedEPSF{lddromadaire.eps scaled 300}}
%\endinsert
\vskip2mm
\Thc Corollaire C {\rm\small (Lemme d'\'elimination de Morse)}|
Soit
${\cal F\/}\!=\!(f, Z)$
une donn\'ee de Morse,
$b, c\!\in\!{\mathop{\Sigma}\nolimits}(f)$
et
$l$
une
$Z$%
-orbite avec
 $\overline{l}\!=\!l\cup\{b, c\}$%
.
On

suppose
$\bigl({\cal P\/}{_{f(c)}}b\bigr)$
et
$$
l={^{^\nwarrow}}c(Z)\cap
b_{_{\searrow}}\!(Z),\
\hbox{\rm  l'intersection \'etant transverse en}\ l \leqno ({\cal I\/}_{{_{d}}l^{c}})$$

Alors pour tout voisinage
$({\cal B\/}\!, {\cal L\/})$
dans
$V$
de la paire
$(\{c\}\cup{_{f(c)}}{b}_{_{\searrow}}\!\!(Z), \overline{l})$%
,\break
il y a un chemin
${\cal F\/}_t\!=\!(f_t, Z_t), t\in{\Bbb R\/}$
 issu de
${\cal F\/}$
de donn\'ees avec~:

$(1)$
les supports des chemins
$f_t$
et
$Z_t$
sont dans
${\cal B\/}$
et
${\cal L\/}$
respectivement.

$(2)$
la restriction
${{\cal F\/}_t}_{|V\!\setminus\!\overline{\cal L\/}}$
est une
$Z_{|V\!\setminus\!\overline{\cal L\/}}$%
-modification.

$(3)$
il y a
$t_0\in{\Bbb R\/}$
tel que si
$t\ne t_0$%
, la donn\'ee
${\cal F\/}_t$
est de Morse. De plus 
$t_{0}$
est point de  travers\'ee de la strate de codimension
$1$
des singularit\'es cubiques%
\footnote{\small
{\it c. a d.\/} ayant le mod\`ele
$k_{s}(x)\!=\!-x_{1}^{2}\!\cdots\!-x_{i}^{2}+x_{i+1}^{2}\!\cdots\!+x_{n}^{2}+x_{n}^{3}+3(s-t_{0})x_{n},
s\in]t_{0}-\epsilon, t_{0}+\epsilon[$%
.}.

$(4)$
si
$t>t_0$%
, l'ensemble des points critiques de
$f_t$
est
${\mathop{\Sigma}\nolimits}(f_t)={\mathop{\Sigma}\nolimits}(f)\!\setminus\!\{b, c\}$%
.

\finc

\npage
%\input D2.tex
%\titrecourant={\hfill\\small R\'eduction du lemme d'\'elimination aux lemmes A, B, C, D\hfill}

%
%\null\vskip3mm
%
  \fancyhead[RE]{Réduction du lemme d'\'elimination aux lemmes A, B, C, D}
  \fancyhead[LE]{\thepage}
\centerline{\bf   \S2 R\'eduction du lemme d'\'elimination  aux lemmes
{\bf A\/}, {\bf B\/}, {\bf C\/}, {\bf D\/}}

Si
$B\!=\!{b}_{_{\searrow}}\!\!(Z), M\!=\!B\cap{\cal U\/}$%
, le {\petcap Compl\'ement C\/} pour la param\'erisation de Morse en
$c$
donne
$(B^{+}, l^{+})$%
, paire de vari\'et\'es contenant
$({_{f(c)}}B, l)$
comme sous-paire ouverte, avec
$\{b, c\}\!\subset\!l^{+}\!\subset\!B^{+}$%
uniques points critiques de
$(f_{|B^{+}}, f_{|l^{+}})$%
.%\hfill\break

On identifie
$l^{+}$%
, le nom des points conserv\'e, \`a un intervalle r\'eel
$I\ni b\!<\!c$%
.

Le {\petcap Compl\'ement B\/} donne alors des vaux de Morse
${\cal \nu\/}_{B^{+}}\!\subset\!{\cal B\/}$
et
${\cal \nu\/}_{l^{+}}\!\subset\!{\cal L\/}\cap B^{+}$
pour les fonctions
$f$
et
$-f_{|B^{+}}$%
, quitte \`a restreindre
$(B^+, l^+)$%
, il y a
$c\!<\!n\in l^+$
avec
$$h_{f, B^+}\!>\!3f(n)-2f(c)\!=\!f(c)+3\epsilon$$

Pour
$\frac{2f(c)+f(n)}{3}\!=\!u'\!<\!\!u\!\!=\!\frac{f(c)+2f(n)}{3}$
et
${\cal N\/}$
disjoint de
$_{u}B^{+}$
et voisinage de
$[c, n]$%
, le {\petcap Corollaire A'\/} donne un chemin
$\tilde{f}_s$
de
$Z$%
-modifications \`a support, dans
${\cal B\/}$
et disjoint de
${\cal N\/}$%
, de
$\tilde{f}_{-1}\!=\!f$
\`a
$\tilde{f}_{\infty}\!=\!\tilde{f}_0$
avec
$$\tilde{f}_0(b)\!<\!f(n)\!=\!\tilde{f}_0(n)\!<\!h_{\tilde{f}_s, B^+}\!=\!h_{f, B^+}$$

Soit
$k_s:l^+\rightarrow{\bf R\/}$
produit par le {\petcap Lemme D\/}  appliqu\'e \`a la restriction
$k\!=\!{\tilde{f_0}}_{|{l}^{+}}$%
. 
 
Le {\petcap Lemme A\/} appliqu\'e \`a
$V\!=\!{\cal L}, f\!=\!-\,{\tilde{f_0}}_{|{\cal L\/}}, e_s\!=\!-k_s$
donne alors un chemin de fonctions
${\tilde{\tilde{f}}_s} : V\rightarrow{\bf R\/}$%
, avec
${{\tilde{\tilde{f}}_s}}_{|V\setminus{\cal L\/}}\!=\!{{{\tilde{f}}_0}\!}_{|V\setminus{\cal L\/}}$%
, de
${\tilde{\tilde{f}}_0}\!=\!\tilde{f}_0$
\`a
${\tilde{\tilde{f}}_{\infty}}\!=\!{\tilde{\tilde{f}}_{1}}$
tel que le chemin
${\cal F}_{t}=(f_{t}, Z_{t})$%
, o\`u~:
$$f_t : V\rightarrow{\bf R\/},\
\hbox{\rm si\/}\
t\!\leq\!0,\ f_t\!=\!\tilde{f}_t,\
\hbox{\rm si\/}\
t\!\geq\!0,\ f_t\!=\!{\tilde{\tilde{f}_t}}$$
et
$Z_{t}$
est obtenu en interpolant%
\footnote{\small
avec
$r\!=\!r_W$
(Cf. {\bf D\'efinion 2\/}
${\romannumeral 1}2)$%
) 
et la fonction
$\beta$
(Cf. le d\'ebut des Appendices), prendre
$Z_t(x)\!=\!(1-\beta(\frac{f_t(x)}{r}))\cdot{\mathop{\rm grad}\nolimits}_{\mu}(f_t)+%
\beta(\frac{f_t(x)}{r})\cdot Z(x)$%
.
}
entre
$V\setminus{{\cal L\/}^{+}}$
et le tube ferm\'e
$T_{{\cal L\/}^{+}}$
de
${\cal L\/}^{+}$
entre
$Z$
et le gradient
${\mathop{\rm grad}\nolimits}_{\mu}(f_t)$
de
$f_t$
pour la m\'etrique
$\mu$
adapt\'ee \`a
$f$%
, convient.
\hfill\findem
\vskip2mm
\TrimTop{-10pct}
\TrimBottom{-7pct}
\centerline{\BoxedEPSF{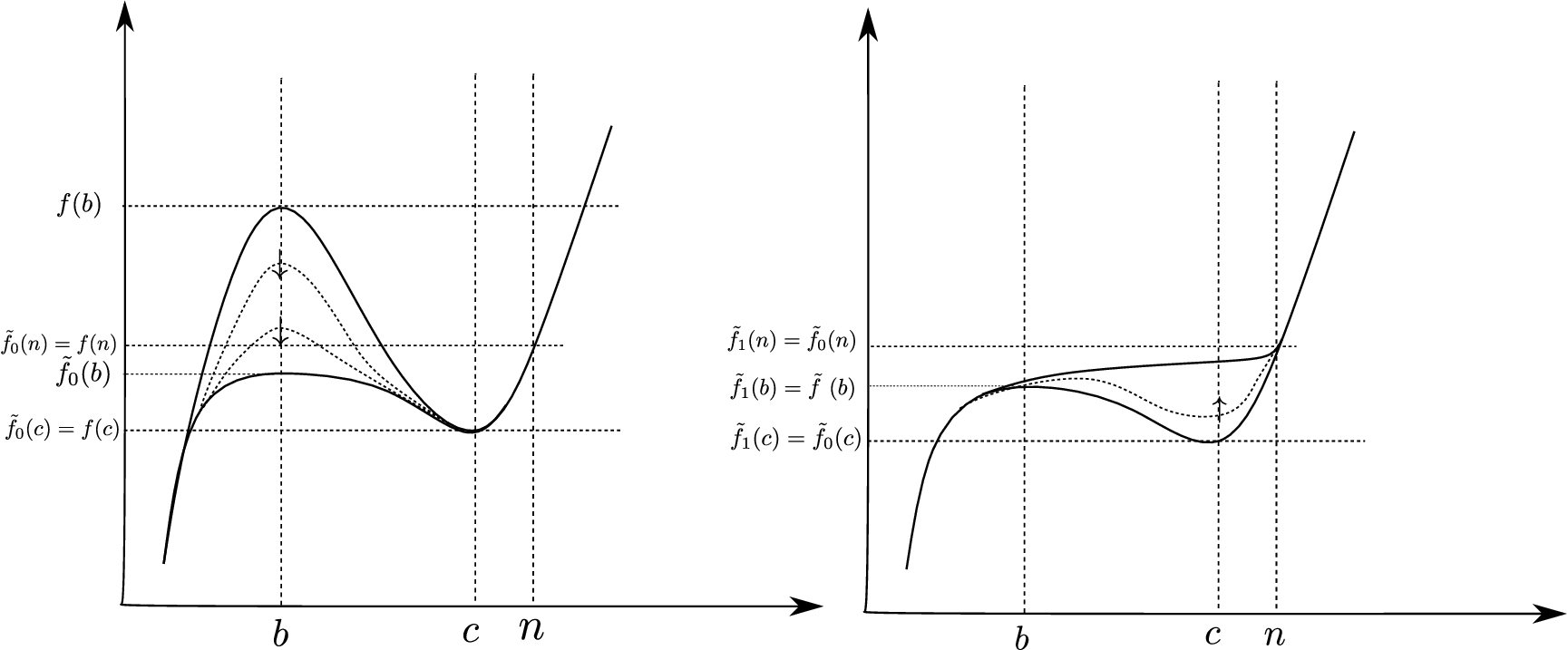 scaled 450}}
%\Lafigps{extension.eps}{500}{extension de nappes}
%\end{figure}
\vskip4mm
\centerline{Baisser la bosse sous le niveau de la nuque\quad\quad\quad puis, appliquer le memme du dromadaire.}
\npage
%\input D3.tex
%\titrecourant={\hfill\\small Commentaires bibliographiques\hfill}
\setcounter{footnote}{0}
  \fancyhead[RE]{\S 3 Commentairs bibliographiques}
  \fancyhead[LE]{\thepage}
%\null\vskip-5mm
\centerline{\bf \S3 Commentaires bibliographiques}
\vskip2mm
En 1925 ({\bf [5]\/}%
$(1)$%
) Morse \'etabli son \og Lemme\fg de forme canonique d'une fonction
$f(x_1,\ldots, x_n)$
de classe
$C^{3}$
en
$0$%
, point critique non d\'eg\'en\'er\'e,
en montrant que
$f$
est translat\'ee d'une \og forme quadratique \`a coefficients variables diff\'erentiables\fg
{\small
Apr\`es factorisation de la premi\`ere variable, la r\'egularit\'e du quotient
par formule de Taylor avec reste ponctuel\footnote%
{\small
A partir de {\bf [113]\/}, pour la forme canonique, il renvoit \`a  {\bf [M]\/}, donnant la premi\`ere
r\'eduction par la formule de Taylor avec reste int\'egral
(La r\'egularit\'e est alors claire et cette premi\`ere \'etape passe \`a param\`etres,
ce n'est pas le cas de la seconde qui reste par diagonalisation de Lagrange.)
}, la premi\`ere r\'eduction s'obtient par r\'ecurrence sur le nombre de variables
\vskip-3mm
$$f(\underline{x})\!=\!f(0, \underline{x'})+a(x_1, \underline{x'})x_1\!=\!%
f(0, \underline{x'})+{{\partial f}\over{\partial x_1}}(0,\underline{x'})x_1+b_{1, 1}(x_1, \underline{x'})(x_1)^{2}\!=\!$$%

$$\!=\!f(0)+\sum_{i,j=2}^nb_{i, j}(\underline{x'})x_i\cdot x_j+\sum_{i=2}^nb_{i, 1}(\underline{x'})x_i\cdot x_1+b_{1, 1}(x_1, \underline{x'})(x_1)^{2}$$%
}
puis  diagonalisant cette derni\`ere par le proc\'ed\'e de Lagrange\footnote%
{\small
Si
$f(\underline{x})\!=\!\sum_{i, j=1}^n b_{i, j}(\underline{x})\cdot x_i\cdot x_j$
avec
$b_{i, j}\!=\!b_{j, i}$
et
$\sum_{i, j=1}^nb_{i, j}(0)\cdot x_i\cdot x_j$
non d\'eg\'en\'er\'ee, apr\`es un \'eventuel changement de variable lin\'eaire pour que
$b_{1, 1}(0)\!\ne\!0$%
, v\'erifier que
$\underline{x}\mapsto\underline{x_2}$
o\`u
$x_{2, 1}\!=\!{{1}\over{\sqrt{|b_{1, 1}|}}}\!\cdot\!x_1, x_{2, i}\!=\!%
-{{b_{1, i}}\over{b_{1, 1}}}\!\cdot\!x_{1}+x_i, 1\!<\!i\!\leq\!n$
est diff\'eomorphisme au voisinage de
$0$%
.\hfill\break
Comme
$f(\underline{x_2})\!=\!\pm (x_{2, 1})^2\!+\!\sum_{i, j=2}^n\!b_{2, i, j}(\underline{x}_2)\cdot x_{2, i}\!\cdot\!x_{2, j}$
et
$\sum_{i, j=2}^n\!b_{2, i, j}(0)\cdot x_{i, 2}\!\cdot\!x_{j, 2}$
est non d\'eg\'en\'er\'ee, le proc\'ed\'e s'it\`ere  pour obtenir la diagonalisation
$f(\underline{x_n})\!=\!\pm (x_{1, n})^2+\cdots+\pm (x_{n, n})^2$%
.
}.

%\vskip-5mm
Si les flots de gradient et les bases des d\'ecompositions en anses
apparaissent dans ce premier article, la g\'eom\'etrie des fonctions de Morse sera en veille\footnote%
{\small
dans l'intervalle le principal domaine d'application des fonctions de Morse sera le calcul variationel abstrait,
notons l'expos\'e {\bf [ST]\/} quasi simultan\'e \`a {\bf [M]\/} de cette \og th\'eorie de de Morse\fg
qui n'utilise pas le lemme de forme canonique, se contentant d'un d\'eveloppement limit\'e.
}
jusqu'\`a  {\bf [113]\/}
en 1960
o\`u il construit sur toute vari\'et\'e compacte connexe, une fonction de Morse {\it polaire\/}\footnote%
{\small
{\it c.ad.\/} avec un seul minimum local et un seul maximum local.
}%
: Quand le sous-niveau
$\{f(x)\!\leq\!c\}$
d'une fonction de Morse
$f$
\`a valeurs critiques distinctes perd une composante connexe quand
$c$
croit passant par la valeur d'un point critique
$A_0$
d'indice un, dit {\it s\'eparant\/},  il exhibe un  minimum
$A_1$
et une d\'eformation \`a support dans un voisinage d'un arc contenant
$\{A_0, A_1\}$
et disjoint de
$\Sigma(F)\setminus\{A_0, A_1\}$
de
$f\!=\!f_0$
\`a
$f_1$
avec
$\Sigma(f_1)\!=\!\Sigma(f)\setminus\{A_0, A_1\}$%
.

Il y lance les bases des lemmes g\'eom\'etriques que Smale utilisera pour le h-cobordisme :
existence de fonctions de Morse \`a valeurs critiques distinctes, modification des gradients
pour que les liaisons ne soient que  vers un point d'indice sup\'erieur\footnote%
{\small
se limitant \`a obtenir que les deux lignes de gradient tendant vers
$A_0$
viennent de deux minimas
$A_1, A_{-1}$
avec
$f(A_{-1})\!<\!f(A_1)$%
, d'o\`u un arc ouvert direct
$\gamma$
passant par
$A_{-1}, A_0, A_1$%
.
}
et, sans lui donner encore de nom,
sa notion de
$f$%
{\it -fonction de base\/}
$\beta$
qu'il ne cessera de d'utiliser dans {\bf [139]\/}, {\bf [117]\/}, {\bf [138]\/}, {\bf [140]\/}, {\bf [150]\/}~:

C'est une fonction
$f: M\rightarrow{\bf R\/}$
sur une vari\'et\'e Riemannien\-ne
$M$
et une sous-vari\'et\'e\footnote%
{\small
localement ferm\'ee, non n\'ecessairement ferm\'ee.
}
$\beta\subset M$
avec
$f$
croissante au sens large sur les germes de g\'eod\'esiques normales\footnote%
{\small
On note
$\epsilon\!>0$
tel que ces germes sont de longueur
$2\epsilon$%
.}
\`a
$\beta$%
, germes dont l'union des images est un tube
${\cal A\/}$
autour de
$\beta$%
. Ainsi~:

{\sl Si
$g :\beta\rightarrow{\bf R\/},\ \theta : {\bf R\/}\rightarrow{\bf R\/}$
sont des fonctions lisses v\'erifiant
$g\leq f_{|\beta}$
avec \'egalit\'e hors d'un compact,
$\theta$
d\'ecroit de
$\theta\!=\!1$
sur
$]\!-\infty, 0]$
\`a
$\theta\!=\!0$
sur
$[\epsilon, +\infty [$
avec
$\theta'\!<\!0$
sur
$]0, \epsilon[$%
, la fonction
${\cal F\/} : M\rightarrow{\bf R\/}$
d\'efinie par
${\cal F\/}\!=\!f$
hors de
${\cal A\/}$
et sur tout rayon normal \`a
$\beta$%
, param\'etr\'e par longeur d'arc de
$r(0)\!=\!b\!\in\!\beta$
\`a
$r\!=\!r(s)\!\in\!{\cal A\/}$%
~:
$${\cal F\/}(r)\!=\!f(r)-\theta(s^2)(f(b)-g(b))\!\leqno\!(\ast)$$
a, si
$\Sigma(f)\!\cap\!{\cal A\/}\!\subset\!\beta$%
, pour uniques points critiques sur
${\cal A\/},\,
\Sigma({\cal F\/}_{|{\cal A\/}})\!=\!\Sigma(g)$
ceux de
$g$%
.
}

La tr\`es faible condition de croissance au sens large sur les rayons suffit
pour la  magie de
$(\ast)$
et permet \`a M. Morse de choisir, hors de voisinages de
$\Sigma(g)$%
, le tube
${\cal A\/}$
fibr\'e par niveaux de
$f$%
. Il le paye par des constructions tr\`es techniques de recollement dans des param\'etrisations de Morse
des points de
$\Sigma(g)$%
.

D'autres notions r\'ecurrentes dans ces articles de Morse sont celles de {\it boule (bowl) (et d\^ome) ascendants et descendant\/}
d'un point critique
$c$
non d\'eg\'en\'e muni d'une param\'etrisation de Morse (et d'un pseudo-gradient
$Z$%
), ce sont~:
$${\cal B\/}_{-}(c)\!=\!{c_{_\searrow}}\!(f, Z),\quad \Bigl(B_{-}(c)\!=\!{\cal B\/}_{-}(c)\cap f^{-1}([f(c'), f(c)])\Bigr)$$
$${\cal B\/}_{+}(c)\!=\!{^{^\nwarrow}c}(f, Z),\quad \Bigl(B_{+}(c)\!=\!{\cal B\/}_{+}(c)\cap f^{-1}([f(c), f(c'')])\Bigr)$$
les nappes et conappes de
$c$
par
$Z$
(tronqu\'ees aux valeurs
$f(c'), f(c'')$
du\footnote%
{\small
sans param\`etres en vue, Morse suppose ses fonctions \`a valeurs critiques distinctes.
}
point cri\-ti\-que
%$c'$
le plus haut et
%$c''$
le plus bas adh\'erent)~:
$\bigl(\overline{B_{-}(c)}, \overline{B_{+}(c)}\bigr)\cap\Sigma(f)\!=\!(\{c', c\}, \{c, c''\})$%
.

{\sl Pour tout voisinage
$N$
de
$B_{\pm}(c)$
et 
$\eta\!\in\!f(B_{\pm}(c))$
il y a 
$\dot{f}$
ayant m\^emes points critiques que
$f$%
, \'egale \`a
$f$
hors de
$N$
et telle que
$\dot{f}(c)\!=\!\eta$%
.}
({\bf [138]\/}, pr\'ecis\'e par {\bf [150]\/}).

La condition\footnote%
{\small
substituant la s\'eparation dans {\bf [113]\/} pour \'eliminer un point d'indice $1$
(avec un minimum).%
}
d'\'elimination d'une paire
$\{c, d\}, f(d)\!<\!f(c)$
de points critiques, introduite dans l'esquisse {\bf [139]\/}\footnote%
{\small
d\'emontrations compl\`etes  dans {\bf [117]\/} qui utilise {\bf [136]\/}, {\bf [140]\/} et {\bf [143]\/}
(annonc\'e dans {\bf [137]\/}).%
}%
, est que {\sl les d\^omes
$B_{-}(c)\cap B_{+}(d)\!=\!\gamma$
se coupent transversalement,  l'intersection \'etant r\'eduite \`a une orbite
$\gamma$
du pseudo-gradient
$Z$%
.%
}

%{\small
Il produit d'abord, sur
$\beta\!=\!B_{\pm}(c)$
(resp. une sous-vari\'et\'e
$\beta\supset B_-(c)\cup\{d\}$
contenant
$d$
et dans laquelle le d\^ome  de
$c$
est ouvert), une fonction
$g$
avec
$c$
unique point critique et
$g(c)\!=\!\eta$
(resp.
$g$
non singuli\`ere) et
$g\!\leq\!f_{|\beta}$
avec \'egalit\'e hors d'un compact, puis l'\'etend par
$(\ast)$%
, une fois prouv\'e que
$\beta$
est base de
$f$%
-fonction.

N'ayant remarqu\'e le {\bf Lemme C\/} la seule construction de
$\beta$
demande une modification du pseudo-gradient  n\'ecessitant de longues pr\'eparations, qui sont reprises dans presque tous les articles,
 pour que le d\^ome de
$c$
entre dans un voisinage de Morse de
$d$
par un  espace de coordon\'ees.

Outre la technicit\'e de la construction {\bf [143]\/} de
$g$
pour l'\'elimination, une grande difficult\'e pour lire  Morse est que, bien que citant {\bf [St]\/}
paru plus de dix ans auparavant, dans presque chaque article Morse\footnote%
{\small
loin de l'aisance de Thom {\bf [T1]\/} pour ces notions dans son article fondateur dix ans avant.%
}
nous demande d'entrer dans ses d\'efinitions
de vari\'et\'es diff\'erentiables, fibr\'es et voisinages tubulaire g\'eod\'esiques, obligeant une lectrice qui sauterait ces \og pr\'ecisions\fg
\`a perdre nombres \og notations permanentes\fg n\'ecessaires  pour entrer dans les parties techniques.
\vskip3mm

\centerline{\bf 1960-1968\/}
\centerline{Le h-cobordisme de Smale et ses expositions par Milnor et Cerf}

Si Morse ne consid\`ere les pseudo-gradients que comme outils de construction,
d\`es 1949 Thom {\bf [T0]\/} les r\'ev\`elent comme apportant une information globale sur la vari\'et\'e
et, en pr\'eliminaire au h-cobordisme, Smale {\bf [Sm0] \/} \'etablit que le gradient d'une fonction de Morse sur une vari\'et\'e
$V$
s'approxime
par un champ de vecteur \`a m\^eme singularit\'es, dont les vari\'et\'es stables et instables
s'intersectent transversalement et dont les ensembles limites sont des singularit\'es.
Il montre qu'un tel champ
$X$
est pseudo-gradient d'une fonction de Morse en construisant inductivement des hypersurfaces
$B_k, 0\!\leq k\!<\!n=\!{\mathop{\rm dim}\nolimits}(V)$
transverses au champ, s\'eparant donc
$V$
en
$G_h$
et
$H_k$
telles que
$G_{k-1}\subset G_{k  }, k\!=\!1,\ldots, n-1$
et
$G_k$
contient les points critiques d'indice
$\leq k$%
, la vari\'et\'e
$G_k$
\'etant obtenue en ajoutant \`a
$G_{k-1}$
des voisinages disjoints des nappes des points critiques d'indice
$k$%
.

Ces voisinage de nappes, diff\'eomorphes \`a des disques, sont des {\it anses\/},
pierre angulaire de la th\'eorie de Smale. Le {\petcap handlebody theorem\/}%
\footnote{\small
{\bf theorem C\/} de {\bf [Sm1]\/}, {\bf 2.1\/} de {\bf [Sm3]\/},
{\bf Lemme fondamental\/} de {\bf [C1]\/} et  {\bf [C2]\/}.
} d'\'elimination de deux anses
$h_{i}, h_{i+1}$
cons\'ecutives%
\footnote{\small
$h_{i}$
et
$h_{i+1}$
attach\'ees \`a
$W_{i-1}$
et
$W_{i}\!=\!W_{i-1}\cup h_{i}$
respectivement,
$W_{i+1}\!=\! W_{i}\cup h_{i+1}$%
.}
{\it en bonne position\/} suit de ce que l'union d'une vari\'et\'e \`a bord
$W$
et d'un disque de m\^{e}me dimension sur des disques de codimension z\'ero de leurs bords est
diff\'eomorphe \`a
$W$%
\footnote{\small
donnant d'abord que
$h_{i}\!\cup\!h_{i+1})$
est un disque, puis
$W_{i+1}\!=\!(W_{i-1}\!\cup\!h_{i})\!\cup\!h_{i+1}\!=\!W_{i-1}\!\cup\!(h_{i}\!\cup\!h_{i+1})$
est diff\'eomorphe \`a
$W_{i-1}$%
.
}%
, m\'ethode perdant tout contr\^{o}le dans l'\'elimination de deux poins critiques en bonne position d'une fonction de Morse%
\footnote{\small
outre ce lemme d'{\it union au bord d'un disque\/}
il faut utiliser deux fois de suite la pas tr\`es explicite construction de {\bf [Sm0]\/}
de passage de d\'ecomposition ordonn\'ee en anses \`a fonction ordonn\'ee esquiss\'ee
au paragraphe pr\'ec\'edent.
}.

C'est sans doute ce qui pousse Milnor ({\bf [MSS]\/} et Cerf ({\bf [CG]\/})%
\footnote{\small
dans ces deux r\'edactions en 1965 et 1968 de cours profess\'es   par Milnor et Cerf en 1963 et 1966
et r\'edig\'ees par des auditeurs de ces cours il est difficile,
sans contact avec auteurs et r\'edacteurs de distinguer ce qui vient de Milnor et Cerf
des \og am\'eliorations\fg des r\'edacteurs!
}
\`a n'y jamais mentionner explicitement les anses, Milnor mettant l'accent sur ses
{\it gradient-like vector fields\/}%
\footnote{\small
nos {\it pseudo-gradients\/},
\'eliminant des singularit\'es en bonne position (Theorem 5.4 de {\bf [MSS]\/})
d'un  pseudo-gradient,
il le modifie au voisinage de la liaison,
sans contr\^{o}ler la fonction.
}
alors qu'ils n'apparaissent pas chez Cerf qui, quitte \`a changer de m\'etrique riemannienne%
\footnote{\small
{\bf Proposition II 2.4\/} de {\bf [CG]\/}. C'est ce qui nous a sugg\'er\'e
notre notion de {\it m\'etrique adapt\'ee\/}
(pr\'ecisant, en imposant des valeurs propres de la Hessienne de norme 1, celle de IV 1 de {\bf [CG]\/})
permettant notre {\bf lemme C\/}
qui nous \'evite toute  construction technique de
{\it double mod\`ele\/}.
},
anticipant sans doute sa {\it th\'eorie de Cerf\/}
de d\'eformation de familles \`a param\`etres de fonctions 
({\bf [C4], [C5]\/} d\'ebut\'ee en {\bf [C3]\/}), se concentre sur les fonctions
et n'utilise que le gradient.
Son approche s'appuie sur une \'etude fine du mod\`ele de Morse,
le modifiant par diff\'eomorphismes de but et source%
\footnote{\small
que, pour ceux de source, il nomme {\it r\'etractions\/}.
}
respectant les niveaux de la fonction 
et la possibilit\'e de,
au lieu de ne les consid\'erer qu'au voisinage des points critiques,
les \'etendre le long de toute {\it nappe d'un point critique\/}
$c$
{\it jusqu'au niveau\/}
$a$%
\footnote{\small
une sous vari\'et\'e de dimension l'indice de
$c$%
, ayant un maximum non d\'eg\'en\'er\'e en
$c$%
, sur laquelle
$f$
est \`a valeurs dans
$[a, f(c)]$%
, non singuli\`ere hors de
$c$
et propre.
},
d'o\`u  un  lemme de descente%
\footnote{\small (Proposition 1 de III de {\bf [CG]\/})}
de d\'emonstration limpide.

Chez Cerf le r\^ole du groupe produit des diff\'eomorphismes
de source et but est essentiel menant au, malheuresement peu cit\'e,
{\it lemme des chemins \'el\'ementaires\/} \'enonc\'e
et montr\'e dans {\bf [C5]\/} et qui permet (haut de la page 434 de {\bf [C4]\/})
de ramener tout chemin g\'en\'erique de fonctions \`a des chemins mod\`eles bien compris.

Pour le lemme d'\'elimination
il  montre que deux singularit\'es en bonne position ont aussi un mod\`ele%
\footnote{\small
a un niveau entre les deux valeurs critiques, il \og plombe\fg la sph\`ere d'attachement
de la nappe montante du point inf\'erieur \`a celle descendante du points sup\'erieur
cr\'eant son \og double mod\`ele de Morse\fg de ces deux points critiques, puis le \og sature\fg en
suivant les lignes de gradients jusqu'\`a deux niveaux entourant
les deux valeurs critiques.
}.
Ceci le ram\`ene \`a un seul cas%
\footnote{\small
celui de son mod\`ele de naissance
$k_{s}(x)\!=\!-x_{1}^{2}\!\cdots\!-x_{i}^{2}+x_{i+1}^{2}\!\cdots\!+x_{n}^{2}+x_{n}^{3}-3sx_{n},
s\!\in\![-1, 1]$
}
qu'il traite explicitement%
\footnote{\small
apr\`es avoir \og engouffr\'e\fg dans un double mod\`ele des points
$(0, \pm1)$
le support d'une fonction cloche
$\omega$
avec
$\omega(x', t)\!=\!1$
pr\`es de
$\{0\}\times[-1, 1]$%
, il prend
$k_{\omega(1-t)}(x), x\in[0, 1]$%
.}.
\vskip3mm

\centerline{\bf 1946-1982\/}
\centerline{${\goth M\/}{\goth J\/}^{2}$, singularit\'es  diff\'erentiables
et m\'ethode du chemin par Samuel, Moser, Mather,\dots}

Contrairement \`a la d\'efinition usuelle des vari\'et\'es diff\'erentiables,
d\'efinies depuis Whitney  par des atlas de cartes
diff\'erentiablement ou analytiquement compatibles (Cf. I de {\bf [Wh]\/}),
l'approche alg\'ebrique par leur (faisceau de) fonctions
semble \^etre syst\'ematiquementexpos\'ee pour la premi\`ere fois par Chevalley
(chapitre III de {\bf [Ch]\/}). Il la dit
dit inspir\'ee de celle des surfaces de Riemann par Weyl ({\bf [We]\/}).

Le \og Lemme
${\goth M\/}{\goth J\/}^{2}$%
\fg
selon ceux qui l'\'enoncent au d\'ebut des ann\'ees 70%
\footnote{\small
{\bf Lemme fondamental\/} de {\bf [CL]\/}, Theorem 1. en tout d\'ebut du cours au C.I.M.E.
{\bf [P]\/}.\hfill\break Il ne semble c\'el\'ebr\'e que par l'\og\'ecole d'Orsay\fgf,
J. Martinet  en 1974 
a un \'enonc\'e \`a la fois plus et moins fort
((b) de la proposition de II 4. de {\bf [Mar]\/},
la traduction anglaise en 1982 de ce cours \`a la P.U.C. de Rio) avec
${\goth M\/}^{2}{\goth J\/}$
en lieu de
${\goth M\/}{\goth J\/}^{2}$%
, mais supposant
${\goth M\/}^{k}\!\subset\!{\goth M\/}^{2}{\goth J\/}$
et le reste dans
${\goth M\/}^{k}$%
.
}
est issu du s\'eminaire de Thom \`a l'I.H.E.S. sur les travaux de Mather.

Sans la m\'ethode du chemin, mais par un proc\'ed\'e d'approximations successives,
il avait \'et\'e d\'egag\'e%
\footnote{\small
Comme signal\'e aux auteurs de {\bf [CL]\/} par B. Teissier que nous remercions de,
pr\`es de 50 ans apr\`es et
par retour d'un courriel,
nous avoir indiqu\'e la r\'ef\'erence {\bf [Sa]\/}.
}
bien avant par Samuel ({\petcap Lemme\/} 2 de {\bf [Sa]\/})
pour donner un mod\`ele alg\'ebriques aux singularit\'es isol\'ees d'une hypersurface d\'efinie sur l'anneau des s\'eries formelles%
\footnote{\small
Si
${\goth P\/}\!=\!(F)\!\subset\!K[[X_{1},\ldots,X_{n}]]\!:=R$
est l'id\'eal premier principal de l'anneau de s\'eries formelles
$R$
d\'efinissant l'hypersurface et le point
$0$
a une singularit\'e isol\'ee alors une puissance
${\goth M\/}^{d}\!\subset\!{\goth J\/}(F)+(F)$
de l'id\'eal maximal de
$R$
est incluse dans l'id\'eal jacobien de la singularit\'e.\hfill\break
Le lemme
${\goth M\/}{\goth J\/}^{2}$
assure alors que la singularit\'e est \'equivalente \`a celle
de l'hypersurface polynomiale d\'efinie par le tronqu\'e \`a
l'ordre
$2d+1$
de la s\'erie formelle
$F$%
.
}.

La m\'ethode du chemin syst\'ematiquement utilis\'ee par J. Mather
dans ses travaux de 1968 \`a 1970
sur la stabilit\'e des applications
$C^{\infty}$%
, par exemple {\petcap Lemma\/} 1 du \S7 de {\bf [Mat]\/}, semble \^{e}tre d\'egag\'ee
avant en 1965 par J. Moser \`a la fin de {\bf [Mo]\/} pour donner une d\'emonstration
alternative de l'unicit\'e, qu'il venait d'\'etablir directement,
\`a isotopie pr\`es des formes volumes de m\^{e}me
int\'egrale sur une vari\'et\'e compacte, en obtenant par la m\^{e}me occasion
l'unicit\'e \`a isotopie pr\`es des 2-formes symplectiques
dans une classe de cohomologie fix\'ee.
\vskip3mm

\centerline{\bf 2012-2014\/}
\centerline{Deux \og reprises\fg de la g\'eom\'etrie de Morse-Cerf par  Laudenbach}

D\`es le d\'ebut de 2012 Laudenbach se repenche sur la th\'eorie de Morse-Cerf pour
d\'emontrer par \og th\'eorie de Cerf\fg le th\'eor\`eme de Reidemeister-Singer%
\footnote{\small {\sl (unicit\'e \`a stabilisation et isotopie pr\`es des scindements de Heegaard des 3-vari\'et\'es)}
publi\'e en 2014 {\bf [L1]}, voir https://arxiv.org/abs/1202.1130, pour les versions ant\'erieures, v4 contient l'\'enonc\'e du theorem 1.3 sans la coupure le rendant incompr\'ehensible dans la version publi\'ee.
}.

%Il prouve que dans un chemin g\'en\'erique
%entre deux fonctions de Morse ordonn\'ees tous
%les points de naissance peuvent se tirer en d\'ebut du chemin
%\`a des niveaux tels
%que, dans ce d\'ebut, les fonctions restent ordonn\'ees,
%il pousse de m\^{e}me  ceux de mort \`a la fin du chemin, puis%
%\footnote{\small
%bien que la d\'emonstration de Cerf vaille \`a
%param\`etres.}
% met des coordonn\'ees sur
%un voisinage d'une nappe descendante pour prouver un lemme de
%descente  passant \`a param\`etres, lui permettant d'ordonne
%les fonctions du chemin entre ces d\'ebuts et fin ordonn\'es.

Il donne%
\footnote{\small
avec force construction \`a l'aide de pseudo gradient d'un mod\`ele de voisinage d'une nappe.
}
un lemme de d\'ecroissance de valeur critique valable
dans une famille
\`a param\`etres%
\footnote{\small
ses r\'ef\'erences, contrairement \`a celles de {\bf [L2]},
n'incluent pas {\bf [CG]} dont la d\'emonstration
de la proposition 1 de III est bien plus \'el\'egante
et fonctionne tout autant \`a param\`etres.
}
et des \og versions \'el\'ementaires\fg de lemmes \'enonc\'es
et montr\'es par Cerf dans le cas simplement connexe
en grande dimension. Il ne d\'etaille que son
``elementary swallow tail lemma''
disant les autres d\'emonstrations \^{e}tre analogues%
\footnote{\small
suivant les leitmotivs
\og se ramener \`a la dimension un\fgf,
m\'ethode du chemin et lemme
${\goth M\/}{\goth J\/}^{2}$%
.
}.
Mais signale, pour le cancellation lemma, sa note {\bf [L2]}.

La \og version fran\c caise abr\'eg\'ee\fgf%
\footnote{\small
\`a part une \og r\'e\'ecriture de l'histoire\fg affirmant que Smale utilisait
le lemme d'\'elimination de Morse alors que, comme rappel\'e plus haut et dans
{\bf [MSS]\/},
Smale  ne travaillait qu'avec les anses,
les fonction de Morse n'\'etant utilis\'ees qu'au tout d\'ebut pour construire
une d\'ecomposition en anses et, \`a part {\bf [113]\/}(5) en indice 0 pour
la construction de fonctions polaires,
tous les articles de Morse sur l'\'elimination des points critiques sont
bien post\'erieurs aux travaux de Smale.\hfill\break
2) une coquille
$h(0)\!=\!f(q)-\epsilon$
au lieu de
$h(0)\!=\!f(p)-\epsilon$
et 3) le manque de pr\'ecaution dans le choix de la fonction croissante
$h_{1}\!\leq\!h$
et la d\'eformation de son lemme d'abaissement faisant qu'il n'y
a aucune raison que dans son chemin entre la fontion initiale et celle
o\`u deux points critiques ont disparu il n'y ait qu'un seul point critique
non d\'eg\'en\'er\'e et que ce point soit cubique.
}
est assez s\'eduisante: Il construit un arc
$A$
en prolongeant, dans la vari\'et\'e instable du point
critique sup\'erieur
$p$%
, la liaison jusqu\`a atteindre
un niveau inf\'erieur \`a l'autre valeur critique
$f(q)$%
,  modifie la liaison au voisinage de
$q$
pour, tout en passant par
$q$%
, \^{e}tre dans la vari\'et\'e stable de
$q$
et un peu au del\`a de
$q$
et munit un voisinage de cet arc contenant un voisinage du d\^{o}me%
\footnote{\small
 en reprenant les termes de Morse rappel\'es ci-dessus.
}
de
$p$
de coordonn\'ees
$(u, y, z)$
telles que
$f(u, y, z)\!=\!h(u)-|y|^{2}+|z|^{2}$%
. Ceci lui permet de consid\'erer la fonction sur ce voisinage comme la trace
d'une famille de fonctions de Morse param\'etr\'ee par la coordonn\'ee
$u$
de l'arc
$A$%
. Il applique alors son lemme pour  descendre la valeur critique
$f(u)\!=\!h(u)$
jusqu'\`a
$h_{1}(u)$
o\`u la fonction
$h_{1}\!\leq\!h$
est strictement croissante et \'egale \`a
$h$
pr\`es du bord de
$A$%
, faisant par l\`a m\^{e}me disparaitre les deux points critiques.
La version anglaise d\'etaill\'ee par contre, oublie que Cerf sait
simplement substituer aux mod\`eles locaux de Morse
ses {\it nappes descendantes\/}
globales et n'ayant  remarqu\'e
notre Lemme C et son compl\'ement(extension canonique d'une nappe transverse)%
\footnote{\small
n\'ecessitant un choix judicieux de m\'etrique riemannienne.
}
qui donne gratuitement sa vari\'et\'e
$W\!=\!\{z\!=\!0\}$%
, a besoin de fastidieux  recollements.
\vskip15mm
%
%\npg
%\input Ref.tex
%\titrecourant={\hfill\\small R\'ef\'erences
%\hfill}
  \fancyhead[RE]{R\'ef\'erences}
  \fancyhead[LE]{\thepage}
\vfill\eject
\centerline{R\'EF\'ERENCES%
\footnote{\small
{\bf [k]\/}(j)
est l'article {\bf [k]\/},  du volume j, des oeuvres compl\`etes de Morse {\bf [M0]\/}.
}}

\vskip2mm
{\hangindent=1cm\hangafter=1\noindent{\bf   [C1]}\
{\petcap J. Cerf}\pointir
 {\sl La th\'eorie de Smale sur le h-cobordisme des vari\'et\'es},\hfill\break
S\'eminaire Henri Cartan, tome 14 (1961-1962), exp. n$^{\circ}$ 11-13, p. 1-23%
, 
{\oldstyle 1962}.
\par}

{\hangindent=1cm\hangafter=1\noindent{\bf   [C2]}\
{\petcap J. Cerf}\pointir
 {\sl Travaux de Smale sur la structure des vari\'et\'es},\hfill\break
S\'eminaire N. Bourbaki, 1961-1962, exp. n $\,^{\circ}$ 230, p. 113-128%
, 
{\oldstyle 1962}.
\par}

{\hangindent=1cm\hangafter=1\noindent{\bf   [C3]}\
{\petcap J. Cerf}\pointir
 {\sl La nullit\'e de $\pi_{0}(Diff(S^{3}))$ 1 Position du probl\`eme},\hfill\break
S\'eminaire H. Cartan 15, 1962-1963, exp. n $\,^{\circ}$ 20, p. 1-29%
, 
{\oldstyle 1963}.
\par}

{\hangindent=1cm\hangafter=1\noindent{\bf   [C4]}\
{\petcap J. Cerf}\pointir
 {\sl Isotopie et pseudo-isotopie},\hfill\break
Proccedingd of the international congress of mathematicians tome 1 - Moscow,  p. 429-437%
, 
{\oldstyle 1966}.
\par}

{\hangindent=1cm\hangafter=1\noindent{\bf   [C5]}\
{\petcap J. Cerf}\pointir
 {\sl La stratification naturelle des espaces de fonctions diff\'erentiables r\'elles et le th\'eor\`eme de la pseudo isotopie},\hfill\break
Pub. Math. IHES, t. 39,  p. 5-173%
, 
{\oldstyle 1970}.
\par}

{\hangindent=1cm\hangafter=1\noindent{\bf   [CG]}\
{\petcap J. Cerf, A. Gramain}\pointir
 {\sl Le th\'eor\`eme du h-cobordisme (Smale)},\hfill\break
Secr. math. Ec. normale sup., cours profess\'e au printemps 1965 \`a la facult\'e des sciences d'Orsay%
\footnote{\small
Cf.
http://www.maths.ed.ac.uk/$\sim$v1ranick/surgery/cerfgram.pdf
}, 
{\oldstyle 1968}.
\par}

{\hangindent=1cm\hangafter=1\noindent{\bf   [CL]}\
{\petcap A.Chenciner, F. Laudenbach}\pointir
 {\sl Singularit\'ees de codimension 1 et chemins \'el\'ementaires d'\'elimination},\hfill\break
C.R.A.S., 270, 1575-1578,
{\oldstyle 1970}.
\par}

{\hangindent=1cm\hangafter=1\noindent{\bf   [Ch]}\
{\petcap  C. Chevalley}\pointir
 {\sl Theory of Lie groups},\hfill\break
Princeton University press,
{\oldstyle 1946}.
\par}

{\hangindent=1cm\hangafter=1\noindent{\bf   [F-R]}\
{\petcap R. Frisson-Roche}\pointir
{\sl La piste oubli\'ee},\hfill\break
Arthaud,
{\oldstyle 1953}.
\par}

{\hangindent=1cm\hangafter=1\noindent{\bf   [L1]}\
{\petcap F. Laudenbach}$\!\!$\pointir$\!\!\!$%
 {\sl A proof of Reidemeister-Singer's theorem by Cerf's methods},\hfill\break
Ann. fac. sc. Toulouse 23 1 (2014), 197-221
{\oldstyle 2014}.
\par}

{\hangindent=1cm\hangafter=1\noindent{\bf   [L2]}\
{\petcap F. Laudenbach}\pointir
 {\sl A proof of Morse's theorem about the cancellation of critical points},\hfill\break
C.R.A.S., 351, 483-488,
{\oldstyle 2013}.
\par}

%{\hangindent=1cm\hangafter=1\noindent{\bf   [M]}\
%{\petcap  J. Milnor}\pointir
% {\sl Microbundles and differentiable structures},
% Princeton September 1961,\hfill\break
% et Collected papers of  John Milnor IV, John McCleary Editor  173-190  - A.M.S.
% {\oldstyle 2009}
%\par}

{\hangindent=1cm\hangafter=1\noindent{\bf   [Mar]}\
{\petcap  J. Martinet}\pointir
 {\sl Singularities of smoth functions and maps},\hfill\break
London Mathematical Society Lecture Note Series 58,
{\oldstyle 1982}.
\par}

{\hangindent=1cm\hangafter=1\noindent{\bf   [Mat]}\
{\petcap  J. Mather}\pointir
 {\sl Stability of $C^{\infty}$ mappings: II Infinitesimal stability implies stability},\hfill\break
Annals of mathematics 89, 254-291,
{\oldstyle 1969}.
\par}

{\hangindent=1cm\hangafter=1\noindent{\bf   [MSS]}\
{\petcap  J. Milnor}\pointir
 {\sl Lectures on the h-cobordism theorem, Notes by L. Siebenmann and J. Sondow},\hfill\break
Princeton university press,
{\oldstyle 1965}.
\par}

{\hangindent=1cm\hangafter=1\noindent{\bf   [MO]}\
{\petcap  M. Morse}\pointir
 {\sl collected works},\hfill\break
Springer,
{\oldstyle 1981}.
\par}

{\hangindent=1cm\hangafter=1\noindent{\bf   [5]}(1)\pointir
 {\sl Relations between the critical points of  a function of n independent variables},\hfill\break
T.A.M.S., 27, 345-396,
{\oldstyle 1925}.
\par}

{\hangindent=1cm\hangafter=1\noindent{\bf   [113]}(5)\pointir
 {\sl The existence of polar non-degenerated functions on differentiable manifolds},\hfill\break
Annals of mathematics,, 71, 352-383,
{\oldstyle 1960}.
\par}

{\hangindent=1cm\hangafter=1\noindent{\bf   [117]}(5)\pointir
 {\sl The elimination of critical points of a non-degenerated function on a differentiable manifold},\hfill\break
Journal d'analyse math\'ematique, 13, 257-316,
{\oldstyle 1964}.
\par}

{\hangindent=1cm\hangafter=1\noindent{\bf   [136]}(5)\
(avec {\petcap W. Huebsch})\pointir
 {\sl Conditioned differentiable isotopies},\hfill\break
Differential analysis colloquium, Bombay, 1-25,
{\oldstyle 1964}.
\par}

{\hangindent=1cm\hangafter=1\noindent{\bf   [137]}(5)\
(avec {\petcap W. Huebsch})\pointir
 {\sl The bowl theorem and a model non-degenerated function},\hfill\break
P.N.A.S., 51, 49-51,
{\oldstyle 1964}.
\par}

{\hangindent=1cm\hangafter=1\noindent{\bf   [138]}(5)\pointir
 {\sl Bowls, f-fiber bundles, and the alteration of critical values},\hfill\break
Anais da academia Brasileira de ciencas, 36, 245-259,
{\oldstyle 1964}.
\par}

{\hangindent=1cm\hangafter=1\noindent{\bf   [139]}(5)\pointir
 {\sl Bowls of a nondegenerated function on a compapct differentiable manifold},\hfill\break
Differential and combinatorial topology (Symposium in honor of M. Morse) Princeton univerity press, 81-103,
{\oldstyle 1964}.
\par}

{\hangindent=1cm\hangafter=1\noindent{\bf   [140]}(5)\pointir
 {\sl Quadratic forms
$\Theta$
and
$\Theta$%
-fiber-bundles},\hfill\break
Annals of mathematics,, 81, 303-340,
{\oldstyle 1965}.
\par}

{\hangindent=1cm\hangafter=1\noindent{\bf   [142]}(5)\pointir
 {\sl The reduction of a function near a nondegenerated critical point},\hfill\break
P.N.A.S., 54, 1759-1764,
{\oldstyle 1965}.
\par}

{\hangindent=1cm\hangafter=1\noindent{\bf   [143]}(5)\
(avec {\petcap W. Huebsch})\pointir
 {\sl A model non-degenerated function},\hfill\break
Revue Roumaine Math. pures Appl., 27, 691-722,
{\oldstyle 1965}.
\par}

{\hangindent=1cm\hangafter=1\noindent{\bf   [147]}(5)\pointir
 {\sl Non-degenerate functions on abstract differentiable manifold},\hfill\break
Journal d'analyse math\'ematique, 19, 231-272,
{\oldstyle 1967}.
\par}

{\hangindent=1cm\hangafter=1\noindent{\bf   [150]}(5)\pointir
 {\sl Bowls, F-fiber bundles, and the alteration of critical values},\hfill\break
P.N.A.S., 60, 1156-1159,
{\oldstyle 1968}.
\par}

{\hangindent=1cm\hangafter=1\noindent{\bf   [M]}\
{\petcap  M. Morse}\pointir
 {\sl Functional topology and abstract variational calculus},\hfill\break
M\'emorial des sciences math\'ematiques  92, Gauthier-Villars, Paris,
{\oldstyle 1939}.
\par}

{\hangindent=1cm\hangafter=1\noindent{\bf   [MC]}\
{\petcap  M. Morse, S. Cairns}\pointir
 {\sl Critical point theory in global analysis and differential topology}, %\hfill\break
Academic press,
{\oldstyle 1969}.
\par}

{\hangindent=1cm\hangafter=1\noindent{\bf   [Mo]}\
{\petcap  J. Moser}\pointir
 {\sl On the volume element of a manifold},\hfill\break
TAMS 120, 286-294,
{\oldstyle 1965}.
\par}

{\hangindent=1cm\hangafter=1\noindent{\bf   [P]}\
{\petcap  V. Po\'enaru}\pointir
 {\sl Lectures on the singularities of $C^{\infty}$ mappings},\hfill\break
Singularities of analytic spaces (A. Tognoli ed.) C.I.M.E., 83-144
{\oldstyle 1974}.
\par}

{\hangindent=1cm\hangafter=1\noindent{\bf   [Sa]}\
{\petcap  P. Samuel}\pointir
 {\sl Alg\'ebricit\'e de certains points alg\'ebro\"{\i}des},\hfill\break
J. Math. pures et appliqu\'ees t. 35 1, 1-6
{\oldstyle 1956}.
\par}

{\hangindent=1cm\hangafter=1\noindent{\bf   [ST]}\
{\petcap  H. Seifert, W. Threlfall}\pointir
 {\sl Variationsrechnung im Grossen}, %\hfill\break
Teubner
{\oldstyle 1938}.
\par}

{\hangindent=1cm\hangafter=1\noindent{\bf   [Sm0]}\
{\petcap  S. Smale}\pointir
 {\sl On gradient dynamical systems},
Ann. of Math. 74  199-206
{\oldstyle 1961}.
\par}

{\hangindent=1cm\hangafter=1\noindent{\bf   [Sm1]}\
{\petcap  S. Smale}\pointir
 {\sl The generalized Poincar\'e's conjecture in dimension greater than four},
B.A.M.S., 66, 373-375
{\oldstyle 1960}.
\par}

{\hangindent=1cm\hangafter=1\noindent{\bf   [Sm2]}\
{\petcap  S. Smale}\pointir
 {\sl Generalized Poincar\'e's conjecture in dimension greater than four},\hfill\break
Annals of mathematics, 74, 391-405
{\oldstyle 1961}.
\par}

{\hangindent=1cm\hangafter=1\noindent{\bf   [Sm3]}\
{\petcap  S. Smale}\pointir
 {\sl On the structure of manifolds},\hfill\break
American journal of mathematics, 84, 387-399
{\oldstyle 1962}.
\par}

{\hangindent=1cm\hangafter=1\noindent{\bf   [St]}\
{\petcap  N. Steerod}\pointir
 {\sl The topology of fiber bundles},\hfill\break
Princeton mathe matical series 14, Princeton university press
{\oldstyle 1951}.
\par}

{\hangindent=1cm\hangafter=1\noindent{\bf   [T0]}\
{\petcap  R. Thom}\pointir
 {\sl Sur  une  partition  en  cellules  associ\'ee \` a une fonction sur une  vari\'ete},\hfill\break
 C.R. Acad. Sci. Paris, 228, pp. 973-975
{\oldstyle 1949}.
\par}

{\hangindent=1cm\hangafter=1\noindent{\bf   [T1]}\
{\petcap  R. Thom}\pointir
 {\sl Quelques propri\'et\'es globales des vari\'et\'es diff\'erentiables},\hfill\break
Comm. Math. Helv., 28, 17-86
{\oldstyle 1954}.
\par}

{\hangindent=1cm\hangafter=1\noindent{\bf   [T2]}\
{\petcap  R. Thom}\pointir
 {\sl Bifurcation des gradients},\hfill\break
Conf\'erence \`a l'I.H.E.S., d\'ecembre 1969
{\oldstyle 1969}.
\par}
{\hangindent=1cm\hangafter=1\noindent{\bf   [We]}\
{\petcap  H. Weyl}\pointir
 {\sl Die Idee der Riemanschen Fl\H ache}, %\hfill\break
Springer
{\oldstyle 1913}.
\par}

{\hangindent=1cm\hangafter=1\noindent{\bf   [Wh]}\
{\petcap  H. Whitney}\pointir
 {\sl Differentiable manifolds},\hfill\break
Annals of mathematics, 37, 645-680
{\oldstyle 1936}.
\par}
\npage
\null\vskip3mm
%\secno=-1
%\titrecourant={\hfill\small Appendices
%\hfill}
\fancyhead[RE]{{\petcap Appendices\/} A ascenseur de val}
  \fancyhead[LE]{\thepage}
\centerline{\petcap Appendices \/}
%\vskip1mm
%\pageno=6
%{\small
\parc
\centerline{\bf Notations\/}

%\vskip1mm

 Soit
$\rho : {\bf R\/}\rightarrow{\bf R\/}$
une fonction lisse \`a support dans
$[0, 1]$
telle%
\footnote{\small
par exemple d\'efinie par si
$x\in]0, 1[$
par
$\rho(x)={{1}\over{\root 4\of{2\pi}}}\cdot{{\exp(-\cotg^2\pi x)}\over{\sin\pi x}}$
et sinon
$\rho(x)=0$%
.
}
que
$\int\rho^2\!=\!1$%
.
\hfill\break
On note
$\beta :{\bf R\/}\!\!\rightarrow\!\!{\bf R\/},\ \beta(x)\!=\!\int_{-\infty}^{x}\!\rho^2$
donc, pour tout
$a, \epsilon\in{\bf R\/}, \epsilon>0$%
, les fonctions
$$\beta_{\epsilon, a}, \alpha^{a}_{\epsilon} : {\bf R\/}\rightarrow{\bf R\/},\quad
\beta_{\epsilon, a}(x)\!=\!\beta({{x-a}\over{\epsilon}}),\
\alpha^{a}_{\epsilon}(x)\!=\!\beta({{a-x}\over{\epsilon}})$$
sont respectivement croissante
de
$0$
sur
$]\!-\infty, a]$
\`a
$1$
sur 
$[a+\epsilon, +\infty[$
et d\'ecroissante de
$1$
sur
$]\!-\infty, a-\epsilon]$
\`a
$0$
sur
$[a, +\infty[$%
. On note
$\alpha_{\epsilon, a}\!=\!1-\beta_{\epsilon, a},
\beta_{\epsilon,}^{a}\!=\!1-\alpha_{\epsilon}^{a}$%
.
%Soit
%$\alpha\!=\!1-\beta$
%et 
%$\alpha_{\epsilon, a}\!=\!1-\beta_{\epsilon, a},
%\beta^{a}_{\epsilon}\!=\!1-\alpha^{a}_{\epsilon}$%
%.
\hfill\break
Ainsi les flots
$\phi_{\epsilon, a,t}, \psi^{a}_{\epsilon, t}\!:\!{\bf R\/}\!\rightarrow\!{\bf R\/}$
des fonctions 
$\beta_{\epsilon, a}, {\alpha}_{\epsilon}^{a}$
vues comme champ de vecteur sur %la droite r\'eelle
${\bf R\/}$%
:

${\romannumeral 1})$
sont de d\'eriv\'ee en espace
${\mathop{\rm D}\nolimits}\phi_{\epsilon, a,t}, {\mathop{\rm D}\nolimits}\psi_{\epsilon, t}^{a}\!=\!%
{\mathop{\rm O}\nolimits}({{1}\over{\epsilon}})$
et respectivement~:

${\romannumeral 2})$
fixent  les demi droites
$]-\infty, a]$
et
$[a, +\infty[$%
.

${\romannumeral 3})$
sont la translation 
par
$t,\ (x, t)\mapsto x+t$
sur
$\{x+t\!\geq\!a+\epsilon\}$
et
$\{x+t\!\leq\!a-\epsilon\}$%
.
\finc
%
%\npage
%\input DAA.tex
\null\vskip3mm
%\secno=-1
%\titrecourant={\hfill\eightrm A Ascenseur de val
%\hfill}
%\vskip3mm
%\pageno=6
%{\eightpoint
\centerline{{\bf A Ascenseur de val\/}  
}

\vskip10mm

\Demd u Lemme A|
Soit
$(a, \epsilon, t)\!=\!(\overline{a}, \overline{\epsilon}, \overline{t})\circ\pi:\!T\!\rightarrow\!{\bf R\/}$
donn\'e par~:
$$\overline{a}\!=\!h_W,\
\overline{\epsilon}\!=\!%
{{1}\over{3}}(h_W-f_{|W})\!=\!{{1}\over{3}}(\overline{a}-f_{|W}),\
\overline{t}\!=\!e-f_{|W}$$
%1
il suffit de d\'efinir le chemin
$f_{e, s}$
par~: si
$v\!\not\in\!T,\ f_{e, s}(v)\!=\!f(v)$
et si
$v\!\in\!T$%
~:
$$f_{e, s}(v)\!=\!\psi_{\epsilon(v), st(v)}^{a(v)}(f(v)) \leqno(A)$$
{\small
En effet
$a, \epsilon, t$
se factorisant par
$\pi$%
, sur chaque fibre
$\pi_{w}\!=\!\pi^{-1}(w), w\in W$%
, l'application
$f_{e, s}$
est compos\'ee au but de
$f_{|\pi_{w}}$
par un diff\'eomorphisme de
${\bf R\/}$
qui, au voisinage
$\{f(v)<a(v))-\epsilon(v)\}$
de
$w$
dans
$\pi_{w}$
est la translation par
$s[e(\pi(v))-f(\pi(v))]$%
.\hfill\carre
}%
\Demd u Corollaire A|
si
$u'\!<\!a\!=\!{{2u'+u}\over{3}}, \epsilon\!=\!{{u-u'}\over{3}}$%
, la compos\'ee au but de la restriction
$f_{|W}$
par les diff\'eomorphismes
$\phi_{\epsilon, a, \beta(\sigma) (u-\kappa)}$
fait le travail sur
$W$%
.

La continuit\'e en
$(f, e)$
du {\petcap Lemme A\/} pour les compacts
$\ L\!=\!{_{a}}W, K\!=\!{_{\kappa-\epsilon}}W$
et le chemin
$e_{\sigma}=\phi_{\epsilon, a, \beta(\sigma) (u-\kappa)}\circ f_{|W}$
donne le
$f_{\sigma}\!=\!f_{e_{\sigma}, 1}$
cherch\'e.
\hfill\findem
\Demd  u Compl\'ement A|
L'hypoth\`ese $e\!\leq\!f_{W}$
\'etant satisfaite sur tout val ind\'ependamment de sa cr\^ete,
on pourra r\'eduire le val \`a loisir.

Dans le support du chemin
$f_{e, s},\ Z^{-1}(0)\!=\!Z_{|W}^{-1}(0)\!=\!%
{\mathop{\Sigma}\nolimits}(f_{|W})\!\subset\!K$%
, la fin du {\petcap Lemme A\/} et la formule
$(A)$
font le travail sur un voisinage ouvert
${\cal S\/}$
de
${\mathop{\Sigma}\nolimits}\!\subset\!W$%
.

Le compact
$L\!\setminus\!{\cal S\/}\!\subset\!W$
est recouvert par un nombre fini d'ouverts
$W_j$
bases de microparam\'etrisations fibr\'ees
$\phi\!=\!\phi_j\!=\!\vartheta\circ\phi_j: W_j\times{\mathop{\rm B}\nolimits}_r\rightarrow
{\cal V\/}_j\!\subset\!\pi^{-1}(W_j)$
du tube, dans lesquelles 
$Z\circ \phi\!=\!(X, Y),\ (w, y)\!\mapsto\!\bigl(X(w, y), Y(w, y)\bigr)$
avec
$X(w, 0)\!=\!Z(w)$%
.

Le champ
$Z$
\'etant tangent \`a
$W$
on a
$Y\!=\!{\mathop{\rm O}\nolimits}(\parallel y\parallel)$
et
$X(f\circ\phi)_{|W_i\times\{0\}}\!>\!0$
 donc comme, pr\`es de
$0\!\in\!{\mathop{\rm B}\nolimits}_r,\ f\circ\phi(w, y)-f\circ\phi(w, 0)\!=\!q(y)$%
, on a
$Y(f\circ\phi)\!=\!{\mathop{\rm O}\nolimits}(\parallel y\parallel^{2})$%
.

Comme
$3\epsilon\!\geq\!h_W\circ\pi-f\!\geq\parallel y\parallel$%
, la propri\'et\'e
${\romannumeral 1})$
du flot
$\psi^{\epsilon}_{a, t}$
et %la formule
$(A)$
donnent
$Y(f_{e, s}\!\circ\!\psi)(w, y)\!=\!{\mathop{\rm O}\nolimits}({{1}\over{\epsilon}}\parallel\!y\!\parallel^{2})\!=\!%
{\mathop{\rm O}\nolimits}(\parallel\!y\!\parallel)$%
, et pour
$\parallel\!y\!\parallel$
assez petit
$Z(f_{e, s})\!>\!0$%
.

Ainsi, quitte \`a r\'eduire le val, on a
$Z(f_{e, s})_{|V\!\setminus\!{\mathop{\Sigma}\nolimits}(f)}\!>\!0$%
.
\hfill\findem
\npg
\null\vskip15mm
%\secno=-1
%\titrecourant={\hfill\\small B Stabilisation infinit\'esimale et stabilisation \hfill}

%\vskip3mm
%\pageno=6 
%{\small
  \fancyhead[RE]{B Stabilisation infinit\'esimale et stabilisation 1 la m\'ethode du chemin de Moser}
  \fancyhead[LE]{\thepage}

\centerline{{\bf B Stabilisation infinit\'esimale et stabilisation\/}  
}

%{\bf 1 M\'ethode du chemin de Moser~: \/}

\vskip5mm
Dans cet appendice on se fixe une sous-vari\'et\'e
$f$%
-directe
$W\!\subset\!V$
de la vari\'et\'e
$V$
et deux fonctions lisses 
$h, k : V\rightarrow{\bf R\/},\, h_{|W\cup{\cal C\/}}\!=\!k_{|W\cup{\cal C\/}}$
\'egales sur
$W$
et sur un voisinage
${\cal C\/}$
dans
$V$ 
d'un ferm\'e
 $C\!\subset\!W$%
.

On note %la diff\'erence
$d=h-k$%
. On a
$d_{|W}=0$
donc pour tout
$w\!\in\!W, d_{w}\!\in\!{\goth W\/}_w$%
.
\vskip7mm
\centerline{\bf 1 La m\'ethode du chemin de Moser\/}
\vskip4mm

Une isotopie%  
\footnote{\small
{\it c. a d.\/} chemin pour
$-\epsilon\!\leq t\!\leq\!1+\epsilon$%
, issu de
$\phi_0\!=\!{\mathop{\rm Id\/}\nolimits}_V$
de diff\'eo\-mor\-phismes de
$V$%
.
}
$\phi_t : V\!\rightarrow\!V, t\!\in\!{\bf R\/}$
fixe
$W$
et est \`a support dans un voisinage
${\cal V\/}_W$
de
$W$
dans
$V$
et, pr\`es de
$W$
 pour
$0\!\leq\!t\!\!\leq1$%
, v\'erifie la relation
$$(k+td)\!\circ\!\phi_t\!=\!k\ \hbox{\rm donc \/} h\circ\phi_1\!=\!k$$
%\hfill\break
si et seulement son champ
$Z_t\!=\!{{d}\over{dt}}\phi_t$
est nul sur
$W$%
, \`a support dans un voisinage
${\underline{\cal V\/}}_W$
de
$\underline{W}\!=\!W\!\times\![0, 1]$
dans
$\underline{\underline{V}}\!=\!V\!\times\!{\bf R\/}$
et ce champ v\'erifie la {\it relation de Moser\/}~:

$$\hbox{\rm pour tout \/}
(v, t)\
\hbox{\rm pr\`es de \/}
\underline{W},\quad Z_t(k+td)(v)\!=\!-d(v) \leqno(M_W)$$

Cette condition \'etant affine en 
$Z_t$%
, on peut d'une part imposer le support de l'isotopie disjoint de
$C$
et d'autre part il suffit que pr\`es de tout
$p\!=\!(w, t)\!\in\!\underline{W}$
il y ait un germe de champ
$Z_{p}$
nul sur
$\underline{W}$
et v\'erifiant la {\it relation de Moser germifi\'ee\/}~:

$$Z_{p}(k+td)\!=\!-d_{p}\!\in\!{\goth A\/}(\underline{\underline{V}})_p\leqno(m_{p, W})$$
{\small
\Dem Soit
$\psi\!=\!\psi^C$
fonction lisse nulle pr\`es de
$C$
et \'egale \`a 
$1$
pr\`es de
$V\!\setminus\!{\cal C\/}$
et soit
$p\!\in\!{\cal V\/}_p\!\subset\!\underline{\underline{V}}$
voisinage de
$p$
tel que 
$(m_{p, W})$
soit l'\'egalit\'e sur
${\cal V\/}_p$
et       
$(\varphi_p\!)_{p\in \underline{W}}$
partition de l'unit\'e subordonn\'ee \`a ce recouvrement ouvert
$({\cal V\/}_p\!)_{p\in \underline{W}}$
de
$\underline{W}$
et
$Z'_t\!=\!\psi\!\cdot\!\bigl(\sum_{p\in\underline{W}}\varphi_p\!\cdot\!Z_{p, t}\bigr)$%
.

Ce champ d\'ependant du temps v\'erifie alors la relation de Moser
$(M_W)$
et son flot est d\'efini sur un voisinage
$\underline{\cal V\/}'_W$
de
$\underline{W}$%
. Si
$\phi\!=\!\phi_W$
est lisse \`a valeurs dans
$[0, 1]$
\'egale \`a
$1$
pr\`es de
$W$
et \`a support dans
$\underline{\cal V\/}'_W$
le champ
$Z_t\!=\!\phi\cdot Z'_t$
produit l'isotopie cherch\'ee.\hfill\findem
}
\vskip3mm

\vskip3mm
  \fancyhead[RE]{B 2 Un lemme ${\goth M\/}{\goth J\/}^{2}$ relatif}
  \fancyhead[LE]{\thepage}
\centerline{\bf 2 Un lemme
${\goth M\/}{\goth J\/}^2$
relatif\/}

\Thc Lemme {\rm (${\goth M\/}{\goth J\/}^{2}$relatif)\/}|
Si pour tout
$y\!\in\!{\mathop{\Sigma\/}\nolimits}_W(h),\
d\!\in\!{\goth W\/}_y{\goth J\/}_{y}^{2}(h)$%
.

Alors il y a des voisinages 
$C\!\subset\!{\cal C\/}'\!\subset\!{\cal C\/}, W\!\subset\!{\cal V\/}'\!\subset\!{\cal V\/}$
de
$C$
et
$W$
dans
$V$
et une iso\-topie
$\phi_t : V\!\rightarrow\!V$
avec~:

$(1)\quad
{\phi_t}_{|W\cup{\cal C\/}'\cup(V\setminus{\cal V\/}))}\!=\!%
{\mathop{\rm Id\/}\nolimits}_{|W\cup{\cal C\/}'\cup(V\setminus{\cal V\/}))}$%

$(2)$\quad pr\`es de
$\Sigma_{W}(h),\
{\mathop{\rm T\/}\nolimits}{\phi_1}\!=\!%
{\mathop{\rm Id\/}\nolimits}$

$(3)\quad h\circ{\phi_1}_{|{\cal V\/}'}\!=\!k_{|{\cal V\/}'}$
\finc
{\small
\Dem Par la m\'ethode du chemin, on v\'erifie les relations de Moser germifi\'ees~:\hfill\break
${\romannumeral 1})$
Si
$w\!\in\!W\!\setminus\!{\mathop{\Sigma\/}\nolimits}_W$
il y a un champ tangent \`a
$W,\ Z_w\!\in\!{\cal Y\/}_w$
tel que
$Z_w(h)(w)\!\ne\!0$%
.

Comme
$d\!\in\!{\goth W\/}_w,\ Z_w(d)\!\in\!{\goth W\/}\_w\!\subset\!{\goth M\/}_w$
et
$Z_w(d)(w)\!=\!0$%
.
Si
${\cal U\/}_{p}$
est un voisinage de
$p\!=\!(w, t_0)$
dans
$\underline{\underline{V}}$
sur lequel
$|Z_w(h)|\!>\!|(t-1)Z_w(d)|$
et
$\rho$
est lisse \`a support dans
${\cal U\/}_w$
\'egale \`a
$1$
pr\`es de
$w$%
, le champ
$Z_{w, t}$
donn\'e par
${Z_{w, t}}_{|{\cal U\/}_p}\!=\!-\rho{{d}\over{Z_w(k)+tZ_w(d)}}\cdot Z_w$
et nul hors de
${\cal U\/}_{p}$%
, v\'erifie %la relation
$(m_{p, W})$%
.\hfill\carre

\rmc Scholie|
Si
$w\!\in\!\Sigma_{W}$
et
$d\in\!{\goth M\/}_{w}{\goth J\/}_{w}^{2}$
alors pour tout
$t\!\in\!{\bf R\/}$%
, on a~:

$${\goth J\/}_w(k+td)\!=\!{\goth J\/}_w(h+(t-1)d)\!=\!{\goth J\/}_w(h)\!\subset\!{\goth M\/}_w$$

\Dem L'inclusion traduit
$w\!\in\!\Sigma_{W}$%
, la premi\`ere \'egalit\'e est tautologique.

Pour la seconde soit
$Z_1,\ldots, Z_n\!\in\!{\cal X\/}_w$
engendrant le
${\goth A\/}_w$%
-module
${\cal X\/}_w$%
.

Alors, si
$m\in\!{\goth M\/}_w$
et
$1\leq k, l\!\leq n$%
, d'apr\`es la r\`egle de Leibnitz~:
$$Z_i\bigl(m\cdot Z_k(h)\cdot Z_l(h)\bigr)
\!=\!Z_i(m)\cdot Z_k(h)\cdot Z_l(h)+m\cdot Z_i\bigl(Z_k(h)\bigr)\cdot Z_l(h)+%
m\cdot Z_k(h)\cdot Z_i\bigl(Z_l(h)\bigr)\!=\!$$
%\vskip-3mm
$$Z_l(h)Z_i(m)\cdot Z_k(h)+mZ_i\bigl(Z_k(h)\bigr)\cdot Z_l(h)+%
mZ_i\bigl(Z_l(h)\bigr)\cdot Z_k(h)\!=\!%
\sum_{j=1}^{n}m_{i,j}\cdot Z_j(h)$$
o\`u, puisque
$W$
est directe dans
$V$
et
$w\!\in\!\Sigma_{W}$%
, on a
$w\!\in\!\Sigma(h)$
donc
${\goth J\/}_{w}\!\subset\!{\goth M\/}_{w}$%
, d'o\`u
$m_{i, j}\!\in\!{\goth M\/}_w$
ainsi
${^{\mathop{\rm t\/}\nolimits}}[Z_1(k+td),\ldots, Z_n(k+td)]\!=\!%
M{^{\mathop{\rm t\/}\nolimits}}[Z_1(h),\ldots, Z_n(h)]$
avec
$M\!=\!\bigl[{\mathop{\rm Id\/}\nolimits}_{n}+(t-1)(m_{i, j})_{0\leq j, i\leq n}\bigr]$%
, donnant l'inclusion
${\goth J\/}_{w}(k+td)\!\subset\!{\goth J\/}_{w}(h)$%
, l'inclusion inverse suivant de ce que, comme
${\mathop{\rm det\/}\nolimits}
M(w)\!=\!1$%
, cette matrice
$M$
est inversible au voisinage de
$w$%
.\hfill\carre\finc
Ainsi si
$y\!\in\!{\mathop{\Sigma\/}\nolimits}_W,%
k-h\!=\!-d\!\in\!{\goth W\/}_y{\goth J\/}_{y}^{2}(k+td)$
et, pour tout
$p\!=\!(y, t)\!\in\!{\underline{\Sigma_W}}$%
, il y a%
\footnote{\small
en notant
$\underline{{\goth W\/}}_p,\underline{{\goth J\/}}_p,\underline{{\goth M\/}}_p$
les germes en
$p$
de fonctions de
${\goth A\/}(\underline{\underline{V}})$
dont les tranches en
$t'$
sont respectivement nulles sur
$W\!\times\{t'\}$%
, en
$p$
et dans
${\goth J\/}_{y}(h)$%
.
}\break
$f_1,\ldots,f_n\!\in\!\underline{{\goth W\/}}_{p}\underline{{\goth J\/}}_{p}\!\subset\!%
\underline{{\goth W\/}}_p\underline{{\goth M\/}}_p$
avec
$Z_{p, t}\!=\!\sum_{i=0}^nf_iZ_i$
d\'efini sur
 ${\cal V\/}_{p}$%
 , nul sur
$\underline{W}$
et v\'erifiant 
$(m_p)$%
.\carre

Soit donc des
$Z_{p, t}$
v\'erifiant, pour tout
$p\!\in\!\underline{W}$%
, la relation de Moser germifi\'ee  
$(m_{p, W})$
 on peut supposer, si
$w\!\not\in\!{\mathop{\Sigma\/}\nolimits}_W\ {\cal U\/}_{w, t}\cap{\cal V\/}_{p}\!=\!\emptyset$%
, ainsi un
$Z_t$
v\'erifiant
$(M_{W})$
associ\'e \`a ces
$Z_{p, t}$%
\'etant nul \`a l'ordre deux sur
$\underline{{\mathop{\Sigma\/}\nolimits}_W}$
a son flot tangent \`a l'identit\'e sur
${\mathop{\Sigma\/}\nolimits}_W$%
.\hfill\carre\findem
}

\vskip3mm
\fancyhead[RE]{B 3 D\'emonstration du {\petcap Lemme\/} B}
\fancyhead[LE]{\thepage}

\centerline{\bf 3 D\'emonstration du
{\petcap lemme B\/}\/}
\vskip2mm
{\parindent=0pt
{\bf N\'ecessit\'e\/} dans
$(1)$\pointir}
Que
$W$
soit
$f$%
-directe est une condition locale,
on peut donc supposer le tube
$f$%
-plat
$\pi: T\rightarrow W$
trivial, muni d'une param\'etrisation fibr\'ee 
$$\phi\!:\!W\!\times\!{\Bbb R\/}^{r}\!\rightarrow\!V, (w, z)\mapsto\phi(w, z)$$

Tout
$X\!\!\in\!{\cal X\/}$
s'\'ecrit 
$X\!=\!Y+\sum_{i=1}^{r}x_{i}Z_{i}, Y\!\in\!{\cal Y\/},%
x_{i}\!\in\!{\goth A\/},
Z_{i}\!=\!T\phi_{\ast}({{\partial}\over{\partial x_{i}}}),%
$
or pour tout
$w\!\in\!W, Z_{i w}\!\in\!{\mathop {\rm ker\/}}(\pi_{w})$%
, donc
$Z_{i}(f)_{w}\!=\!0$
et
$\rho(X(f))\!=\!\rho(Y(f))\!\in\!{\goth J\/}(f_{|W})$%
.\hfill\carre

{\parindent=0pt
{\bf Suffisance\/} dans
$(1)$\pointir}
$W$
\'etant paracompacte, on peut supposer
$V$
munie, pr\`es de
$W$%
, d'une m\'etrique riemannienne.
Le tube
$f$%
-plat cherch\'e sera g\'eod\'esiquement radial construit \`a partir
d'une section%
\footnote{\small
de l'application quotient
$T_{|W}V\rightarrow T_{|V}/TW$%
.
}
$s: T_{|W}V/TW\rightarrow TV_{|W}$
du fibr\'e normal \`a
$W$
dans
$V$
telle que pour tout
$w\!\in\!W, z\!\in\!T_{w}V/T_{w}W, s(z)(f)(w)\!=\!0$%
.

Ce fibr\'e normal a une section%
\footnote{\small
quelquonque, sans la condition
$s(z)(f)(w)\!=\!0$
ci-dessus, possible
car
$W$
est paracompacte.
}
$s_{\emptyset}$
il suffit d'\'etablir la~:

\rmc Scholie|
Soit
$C, D\!\subset\!W$
un ferm\'e et un compact de
$W$
et
$s_{C}: T_{|W}V/TW\rightarrow TV_{|W}$
une section telle que pour tout
$c'\!\in\!W$
pr\`es de
$C,\, z'\!\in\!T_{c'}V/T_{c'}W$
on a
$s_{C}(z')(f)(c')\!=\!0$%
.

Alors il y a une section
$s_{C\cup D}: T_{|W}V/TW\rightarrow TV_{|W}$
telle que, pr\`es de
$C,\ s_{C\cup D}\!=\!s_{C}$
et pour tout
$c''$
pr\`es de
$C\cup D$
et
$z''\!\in\!T_{c''}V/T_{c''}W$
on a
$s_{C\cup D}(z'')(f)(c'')\!=\!0$%
.

\Dem
Le compact
$D$
est recouvert par un nombre fini de compacts
$D_{k}$
inclus
dans des ouverts
$U_{k}$
de param\'etrisations fibr\'ees
$\phi: U_{k}\!\times\!{\Bbb R\/}^{r}, (u, z)\mapsto \phi(u, z)$
d'un (micro-)tube g\'eod\'esique normal d\'efini par
$s_{C}$%
. Il suffit de traiter le cas o\`u
$D\!=\!D_{k}$
est l'un de ceux-ci.

Comme
$W$
est
$f$%
-directe il y a des champs de vecteurs tangent \`a
$W,\ Y_{i}\!\in\!{\cal Y\/}$
tels que pour
$i\!=\!1,\ldots,r$
on ait, si
$Z_{i}\!=\!\phi_{\ast}{{\partial }\over{\partial z_{i}}}-Y_{i},\, \rho(Z_{i})\!=\!0$
.

Soit
$\alpha, \beta: U\rightarrow[0, 1], \alpha(C)\!=\!1\!=\!\beta(C\cup D)$
lisses \`a support dans des voisinages de
$C\cap U$
et
$(C\cup D)\cap U$
respectivement. La section
$s_{C\cup D}$
d\'efinie en modifiant
$s_{C}$
sur
$U$
par
$s_{C\cup D}({{\partial\phi}\over{\partial z_{i}}})\!=\!%
\beta(\alpha{{\partial \phi}\over{\partial z_{i}}}+(1-\alpha)Z_{i})+
(1-\beta)\phi_{\ast}({{\partial}\over{\partial z_{i}}})$
convient.\hfill\carre
\finc

{\parindent=0pt
{\bf N\'ecessit\'e\/} dans
$(2)$\pointir}
les conditions
$({\cal D\/}), (D_{\goth W\/}), (P)$
\'etant locales%
\footnote{\small
et invariantes par diff\'eomor\-phismes de la source.
}, on peut supposer
$V\!=\!W\!\times\!{\cal N\/}$
o\`u
${\cal N\/}$ est un voisinage de
$0\!\in\!{\Bbb R\/}^{r}$
et, pour une famille
$q_{w}\!:{\Bbb R\/}^{r}\!\rightarrow\!{\Bbb R\/}, w\!\in\!W$
de formes quadratiques d\'efinies positives que~:
$$f : (V, W)\!=\!({\cal W\/}\times({\cal N\/}, \{0\})%
\rightarrow{\bf R\/},\
f(w, z)\!=\!f(w, 0)+q_{w}(z)$$

Les coordonn\'ees
$z_{i}, i\!=\!1,\!\ldots,\!r$
${\goth W\/}$
engendrent%
\footnote{\small
comme
${\goth A\/}$%
-module:
${\goth W\/}=\langle z_{1},\ldots, z_{r}\rangle\!:=%
{\goth A\/}z_{1}+\cdots+{\goth A\/}z_{r}$%
.
}
et
${\cal X\/}$
l'est par
${\cal Y\/}$
et les
$Z_{i}\!=\!%
{{\partial}\over{\partial z_i}}$%
.

Comme
$Z_{i}(f)(w, 0)\!=\!Z_{i}(q_{w})(0)\!=\!0$
la condition
$({\cal D\/})$
est v\'erifi\'ee.\hfill\carre

On a%
\footnote{\small
car
$q_{w}$
est non-d\'eg\'en\'er\'ee.
}
$\langle Z_1(f),\ldots, Z_r(f)\rangle\!=\!%
\langle z_1(f),\ldots, z_r(f)\rangle\!=\!{\goth W\/}$
d'o\`u 
$(D_{\goth W\/})$%
.\hfill\carre

Enfin si
$w\!\in\!{\mathop{\Sigma\/}\nolimits}_{W}$
et si
$X\!=\!Y+Z\!\in{\cal X\/}_{w}\!=\!{\cal Y\/}_{w}+\langle Z_{1},\ldots,Z_{r}\rangle$
est
$f$%
-plat et, comme
$Z$
l'est, on a
$Z(f)\!\in\!{\goth W\/}$%
, donc
$Y(f)\!\in\!{\goth W\/}$
ainsi, par d\'efinition de
${\cal Y\/}, Y(Z(f)), Y(Y(f))\!\in\!{\goth W\/}$
donc, comme les
$Z_{i}$
commutent%
\footnote{\small
car les
$z_{i}$
sont des coordonn\'ees normales \`a
$W$%
.
}
sur
$W$
avec les vecteurs tangents \`a
$W$%
, on a
$X(X\!(f))(w)\!=\!X(Y\!(f)+Z(f))(w)%
\!=\!Y(Y\!(f)(w)+Y(Z(f))(w)+Z(Y\!(f))(w)+Z(Z\!(f))(w)%
\!=\!Y(Y\!(f)(w)+2Y(Z(f))(w)+Z(Z\!(f))(w)$
et
$X(X(f))(w)\!=\!Z(Z(f))(w)\!=\!Z(Z(q_{w}))\geq0$
puisque la forme quadratique
$q_{w}$
est positive d'o\`u la relation
$(P)$%
.\hfill\carre\findem

{\parindent=0pt
{\bf Suffisance\/} dans
$(2)$\pointir}
Soit
$\pi\!:\!T\!\rightarrow\!W$%
, tube
$f$%
-plat
autour de
$W$
donn\'e par
$(1)$%
.

Ainsi tout point
$w\!\in\!W$
est point singulier de la restriction de
$f$
\`a
$\pi^{-1}(w)$%
.\hfill\break
On note
$q_{Hw}$
la Hessienne de cette restriction.

Si
$w_{0}\!\in\!{\mathop{\Sigma}\nolimits}_{W}$
est point singulier de
$f_{|W}$
cette Hessienne
$q_{Hw_{0}}:{\mathop{T}\nolimits}_{w_{0}}\pi^{-1}(w_{0})\rightarrow{\Bbb R\/}$
 est selon
$(D_{\goth W\/})$
non d\'eg\'en\'er\'ee et, d'apr\`es
$(P)$
 positive. Il en est de m\^{e}me pour tout
$w\!\in\!{\cal S\/}_{w_{0}}$
dans un voisinage de
$w_{0}$% 
.

Par convexit\'e des formes quadratiques d\'efinies positives%
\footnote{\small
et puisque
$W$
est  paracompacte.
}, ce germe en
${\mathop{\Sigma}\nolimits}_{W}$
de famille de formes quadratiques d\'efinies positives
$((q_{w})_{w\in\cup{\cal S\/}_{w_{0}}})_{w_{0}\in{\mathop{\Sigma}\nolimits}_{W}}$
s'\'etend en une m\'etrique euclidienne
$\tilde{q}: E\!\rightarrow{\Bbb R\/}$
sur le fibr\'e normal
$\tilde{\pi}\!: E\!=\!{\mathop{\rm T\/}\nolimits}V\big/%
{\mathop{\rm T\/}\nolimits}W\!%
\rightarrow\!W$%
.

On note
${\cal E\/}$
un voisinage de
$W$
dans
$E$
sur lequel est d\'efini un plongement ouvert
$\vartheta:{\cal E\/}\rightarrow V$
d'image un voisinage de
$T$
tel que, sur
$\vartheta^{-1}(T)$
on ait
$\tilde{\pi}\!=\!\pi\circ\vartheta$%
.

Comme
$f$
est, en chaque
$w\!\in\!W$%
, singuli\`ere le long des fibres de
$\pi$%
, la formule de Taylor \`a l'ordre
$2$
donne
$F_2 : X\rightarrow{\cal Q\/}(E)$
lisse \`a valeur dans les formes quadratiques sur
$E$
telle que pour tout
$x\!=\!(w, z)\!\in\!{\cal E\/}\ [w\!=\!\tilde{\pi}(x)\!\in\!W,\ z\!\in\!{\tilde{\pi}}^{-1}(w)]$%
, on a~:
$$f\circ\vartheta(w, z)\!=\!f(w)+F_2(w, z)(z)\!=\!%
f\circ\vartheta\circ{\tilde \pi}(x)+\tilde{q}_{w}(z)+d(w, z)$$

Avec les notations du {\petcap Lemme\/}
${\goth M\/}{\goth J\/}^2$
relatif~:

\Rmc Affirmation |
Le {\petcap Lemme\/}
${\goth M\/}{\goth J\/}^2$
relatif s'applique au reste~:
$$d : X\rightarrow{\bf R\/},\, d(w, z)\!=\!F_2(w, z)(z)-\tilde{q}_w(z)$$
\finc

{\small
\Dem
Soit
$w$
pr\`es de
$w_{0}\!\in\!{\mathop{\Sigma\/}\nolimits}_W$%
. Le choix de la m\'etrique  de
$E$
donne
${\tilde{q}}_{w}\!=\!F_2(w, 0)$%
. Rappelons que
$z\in\!{\goth W\/}$
et, d'apr\`es 
$(D_{\goth W\/}),\, {\goth W\/}_{w_{0}}\!\subset\!{\goth J\/}_{w_{0}}$%
,  
donc%
\footnote{\small
car, par  Taylor \`a l'ordre
$1,\ F_{2}(w, z)-F_{2}(w, 0)\!=\!F_{1}(w, z)(z), F_{1}(w, z)\!\in\!%
{\mathop{\rm End}\nolimits}(E, {\cal Q\/})$
et lisse.
}~:
$$F_2(w, z)(z)-\tilde{q}_w(z)\!=\![F_2(w, z)-F_2(w, 0)](z)\!\in\!%
{\goth W\/}_y({\goth W\/}_{w_{0}}^2)\!\subset\!%
{\goth W\/}_{w_{0}}{\goth J\/}_{w_{0}}^{2}$$

\Demf
Si 
${\cal V\/}'$
et
$\phi_t: {\cal E\/}\rightarrow{\cal E\/}$
sont le voisinage de
$W$
et l'isotopie donn\'ee par le {\petcap Lemme\/}
${\goth M\/}{\goth J\/}^2$
relatif et
$E_{2r'}\!\subset\!{\cal V\/}'$
alors
$\vartheta_1\!=\!{\vartheta\circ\phi_1}_{|E_{2r'}}$
est un plongement de microvoisinage tubulaire de val pour
$f$
cherch\'e
.
\hfill\findem
}

\rmc Remarque|
Si 
$w\!\in\!{\mathop{\Sigma\/}\nolimits}(f_{|W})$
est non d\'eg\'en\'er\'e alors il est
${\goth W\/}$%
-d\'efini positif si et seulement si, comme point critique de
$f$
il est non d\'eg\'en\'er\'e de m\^eme indice que pour
$f_{|W}$%
.
\finc
{\small
\Dem Comme
$w\!\in\!{\mathop{\Sigma\/}\nolimits}(f_{|W})$
est non d\'eg\'en\'er\'e on a
${\goth J\/}_w(f_{|W})\!=\!{\goth M\/}_w(W)$
donc
$(D_{\goth W})$
implique
${\goth J\/}_w(f)\!=\!{\goth M\/}_w(V)$
({\it c. a d.\/}
$w\!\in\!{\mathop{\Sigma\/}\nolimits}(f)$
non d\'eg\'en\'er\'e) et, par
$(P)$%
, si
$N\!\subset\!{\mathop{\rm T\/}\nolimits}_wV$
est sous-espace
$q_w$%
-d\'efini n\'egatif maximal alors
$N\!\subset\!{\mathop{\rm T\/}\nolimits}_wW$
donc
${\mathop{\rm Ind\/}\nolimits}_{w}(f)\!=\!{\mathop{\rm Ind\/}\nolimits}_w({f_{|W}})$%
.\hfill\carre

R\'eciproquement, si
$w\!\in\!{\mathop{\Sigma\/}\nolimits}(f)$
est non d\'eg\'en\'er\'e de m\^eme indice que
$w\!\in\!{\mathop{\Sigma\/}\nolimits}(f_{|W})$%
%\nobreak
alors%
\footnote{\small
car tout
$X_{w}$
dans l'orthogonal pour
$q_{Hw}$
de
${\mathop{\rm T\/}\nolimits}W_{w}$
s'\'etend en un champ
$f$%
-plat et en notant
$q_w^{\perp_b}$
la restriction de
$q_{Hw}$
\`a cet orthogonal.
}
$q_w^{\perp_b}$
est d\'efinie positive (donc non d\'eg\'en\'er\'ee) d'o\`u
$(P)$
et
$(D_{\goth W\/})$
et
${\goth J\/}(f_{|W})_{w}\!=\!{\goth M\/}_{w}$
d'o\`u
$({\cal D\/})$%
.\hfill\carre\findem
}
\null\vskip7mm
%\secno=-1
%\titrecourant={\hfill\\small C Canonicit\'e de l'extension d'une nappe transverse\hfill}
%\vskip3mm
%\pageno=6
%{\small
%\secno=-
%\centerline{\petcap Appendices}
%\vskip1mm
%\pageno=
\fancyhead[RE]{C Canonicit\'e de l'extension d'une nappe transverse}
  \fancyhead[LE]{\thepage}
\centerline{{\bf C Canonicit\'e de l'extension d'une nappe transverse\/}  
}

\vskip5mm
\parc
Avec les notations du {\bf Lemme C\/}, le rayon
$\gamma\cap{\cal U\/}^{\delta}$
a pour voisinage dans
${\cal U\/}^{\delta}$
le c\^{o}ne
$$V^{\delta}(\gamma)\!=\!N^{\delta}+\{\tau\cdot p+r\,;\,r\!\in\!R^{\delta}, s^2\mu(r)\!\leq\!\tau^2\rho^2\}$$
\finc
\Demd   u Lemme C|
$(1)$
suit de la transversalit\'e en
$s\cdot p$
de
$M_s$
et
$P_s$%
, transversalit\'e qui d\'ecoule de celle de
$M$
et
$P$
en
$s\cdot p\!\in\!\gamma\!\subset\!M\cap P$%
.
\hfill\carre
\Rmc {\bf Affirmation\/}|
$({\romannumeral 1})$
Pour la bonne d\'efinition de
$\chi$%
, il suffit de prendre
$0\!<\!\delta\!\leq\!\sqrt{s\nu}$

$({\romannumeral 2})$
On a
$V^{\delta}(\gamma)\!\subset\!\chi(N^{\delta}\times]0, \delta[)\!\subset\!M$
et l'image de
$\chi$
est tangente \`a
$Z^U$%
.
\finc

\preuved e $({\romannumeral 1})$|%
$({{\tau}\over{s}})\delta\!<\!{{\delta^2}\over{s}}\!\leq\!\nu$
ainsi
$\theta$
\'etant lisse,
$\chi$
est d\'efinie et lisse.

Comme
$\chi\!=\!{\mathop{\rm Id}\nolimits}+{\mathop{\rm O}\nolimits}(\tau^2\sqrt{\mu^N(x)})$%
, pour
$\delta$
assez petit,
$\chi$
est un plongement.\hfill\carre

\preuved e $({\romannumeral 2})$|%
Si
$0\!<\!\tau\!<\!\delta,\ t\!=\!\log({{s}\over{\tau}}),\ n\!\in\!N^{\delta},\ r\!\in\!R$
et
$u\!=\!n+\tau\cdot p+r$%
, on a~:
$$\Lambda_{t}(u)\!=\!{{\tau}\over{s}}\cdot n+s\cdot p+{{s}\over{\tau}}\cdot r$$

Si
$u\!=\!\chi(n, \tau),\, u'\!=\!\Lambda_{t}(u)\!=\!%
{{\tau}\over{s}}\!\cdot\!n+s\cdot p+\theta({{\tau}\over{s}}\!\cdot\!n)\!\in\!M$
et
$\Phi^U_{t}(u)$
est  d\'efini et,\break
par invariance de
$M$
par
$\Phi^U_t,\, \chi(n, \tau)\!\in\!M$%
, d'o\`u  la seconde inclusion. D'autre part~:\hfill\break
$s\tau{{d}\over{d\tau}}\chi(n, \tau)_{\tau_0}\!=\!%
s\tau{{d}\over{d\tau}}\Phi^U_{-t}(u')_{\tau_0}\!=\!%
Z^U(\chi(n, \tau_0))$
et
$Z^U$
est tangent \`a l'image de
$\chi$%
.\hfill\carre

Si
$u\!\in\!V^{\delta}(\gamma)\cap M,\ \Lambda_t(u)\!\in\!s\cdot p+(N^{\nu}+R^{\rho})_0$
et
$\Phi^U_t(u)$
est encore d\'efini. Par invariance de
$M$
par
$\Phi^U_t,\ \Lambda_t(u)\!\in\!M$%
, donc d'apr\`es
$(1),\ {{s}\over{\tau}}\cdot r\!=\!\theta({{\tau}\over{s}}\cdot n)$
et
$u\!=\!\chi(n, \tau)\!\in\!\chi(N^{\delta}\times]0, \delta[)$
et la premi\`ere inclusion.\hfill\carre

$(3)$
du {\petcap Lemme C\/}
s'obtient, quitte \`a \'eventuellement r\'eduire
$\delta$%
, en remarquant~:
$$f\circ\psi\circ\chi(n, \tau)\!=\!%
-\mu^N(n)+\mu^P(s\cdot p)+{\mathop{\rm O}\nolimits}\bigl(\tau^4\mu^N(n)\bigr)$$
\hfill\carre\findem

\Demd u Compl\'ement C|%
D'apr\`es la premi\`ere inclusion de
$({\romannumeral 2})$
il n'y a\break
pas de voisinage
$M^{+}\!\!$
de
$\{0\}\!\cup\!M'$
dans
$M_{\chi}$
de topologie de vari\'et\'e, celle induite par l'immersion
$\xi$%
, si et seulement si il y a une suite
$m_n\!\!\in\!M'\!\setminus\!\!V(\gamma)$
qui tende vers
$0$%
.

Si
$c_n(t)$
est la courbe int\'egrale de
$Z^U$
avec
$c_n(0)\!=\!m_n$%
, il y a un
$t_n\!>\!0$
tel que
$c_n(t_n)\!=\!m'_n\!\in\!{\mathop{\rm Fr\/}\nolimits}(N^{\delta}\times P^{\delta})$%
. Cette suite a une valeur d'accumulation
$m'_{\infty}$
qui d'une part, puisque
$m_n\!\in\!M$
et par
$\Phi^U_t$%
-invariance
$m'_n\!\in\!M$
est%
\footnote{\small
puisque
$q(m_n)\!>\!0$
car
$m_n\!\in\!M'$
et
$q$
croit sur les courbes int\'egrale de
$Z^U\!=\!{\mathop{\rm grad\/}\nolimits}_{\mu}(q)$%
.
}
dans
$M'$
et d'autre part, comme
$m_n$
tend vers
$0$
et
$V^{\delta}(\gamma)\subset{\mathop{\rm Im}\nolimits}(\chi)$%
, doit \^etre dans
${\mathop{\rm Fr}\nolimits}(P^{\delta}\setminus V(\gamma))$%
, donc distincte de
$\delta\cdot p$
et
$\bigl({\bf R\/}_{+}^{\ast}\cdot m'_{\infty}\bigr)\cap U\subset M\cap P\setminus\gamma$
est un autre rayon d'intersection.

\hfill\findem
%\secno=-1
%\titrecourant={\hfill\\small B Stabilisation infinit\'esimale et stabilisation %\hfill}

%\vskip3mm
%\pageno=6 
%{\small

%
\npage
%\input DAD.tex
%\titrecourant={\hfill\\small Lemme du dromadaire\hfill}
\fancyhead[RE]{D Lemme du dromadaire}
  \fancyhead[LE]{\thepage}

\centerline{\bf D Lemme du dromadaire}
\vskip3mm
%\pageno=6
%{\eightpoint
%\centerline{{\bf A\/} Ascenseur de val et fonctions dromadaires}

%\vskip1mm

Sur 
$I$
soit huit points
$b'\!<\!b\!<\!c'\!<\!c\!<\!d\!<\!d'\!<\!n\!<\!n',\ k(d')\!<\!k(c)$%
.

Un polyn\^ome cubique
$p$
tel 
$p_i(X)\!=\!2i(X-d)-3(c-d)(X-d)^2+2(X-d)^3,\ $
$i\in\{0, 1\}$%
, est 
{\it rahla\/}%
\footnote{\small
selle de dromadaire des sahariens voir {\bf [FR]\/}.
}
d'une fonction
$k$
({\it en
$c, d, n$%
-dromadaire\/}) v\'erifiant
$(D)$
si~:
$$p'(n)\!>\!0,\ p(d)\!=\!0\quad
\hbox{\rm et \/}\quad
p'^{-1}(0)\in\bigl\{\{c, d\}, \emptyset\bigr\}$$

%{\eightpoint
\Thc Scolie|
Soit une rahla
$q$
de
$k$%
.
Il y a alors
$0\!<\!\eta$
et une fonction
$g_1$
 avec~:

$(1)\quad
{g_1}_{|[c', d']}\!=\!k(c)+\eta+\eta{q}_{|[c', d']}$

$(2)\quad
k\!\leq\!g_1\
\hbox{\rm avec \'egalit\'e sur \/}\
I\setminus]b, n[$

$(3)\quad
{g'_1}_{|I\setminus[c, d]}>0$

$(4)\quad
g_1(c)\!<\!g_1(n)$

\finc
%\midinsert
\TrimTop{-10pct}
\TrimBottom{-5pct}
\TrimLeft{-5pct}
\TrimRight{-5pct}
\centerline{\BoxedEPSF{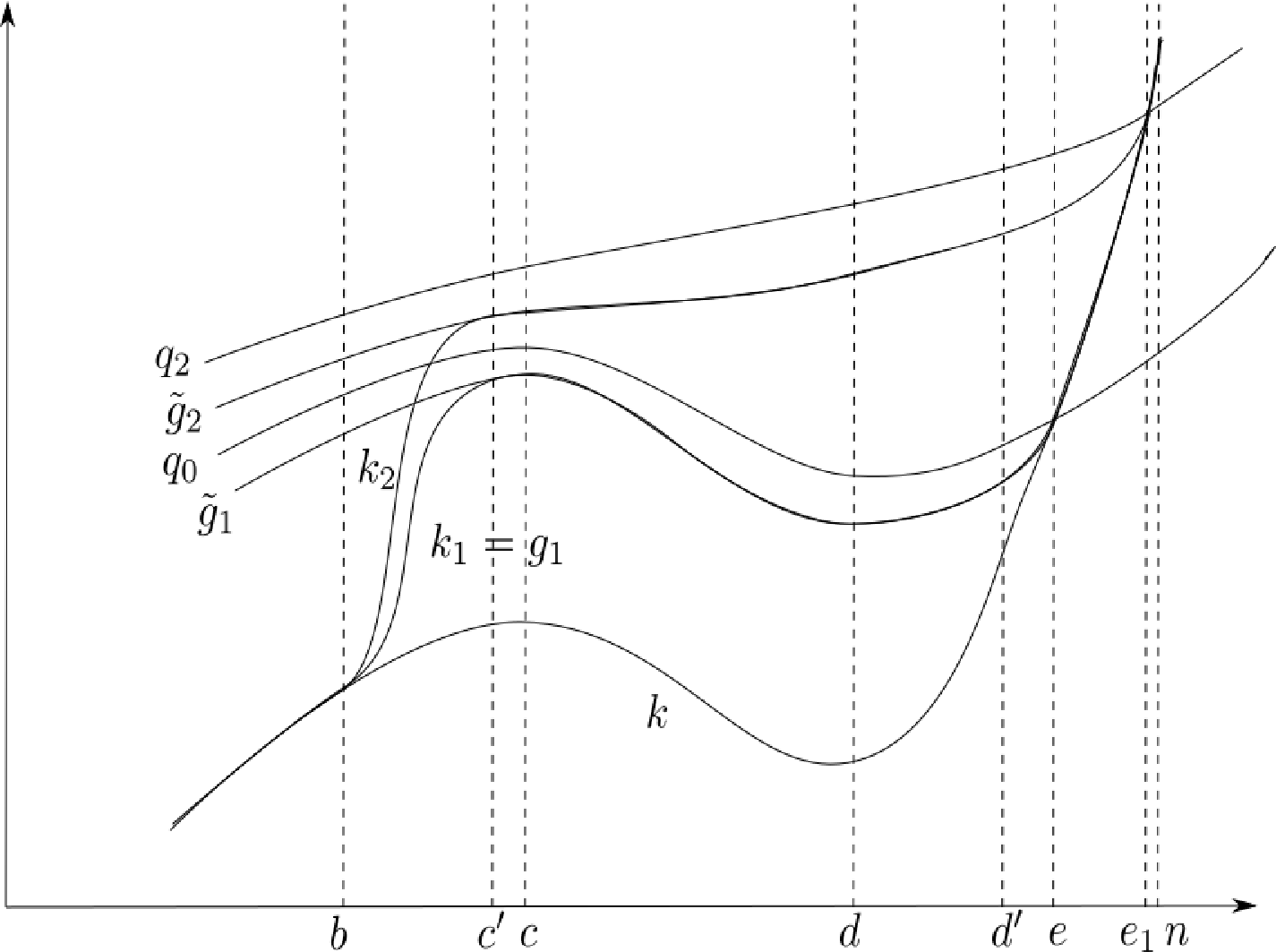 scaled 250}}
%\endinsert
\vskip5mm
\Dem Soit
$q_{1}\!=\!k(c)+\eta(1+q_{|I}):I\rightarrow{\Bbb R\/}$%
. Si
$0\!<\!\eta$
est assez petit
 on a
$k\!<_{[b', d']}{q_1}$
et
$0\!<_{[d', n']}{q_1'}\!<_{[d', n']}{k'}$
et, puisque
$k(c)\!<\!k(n)$
et
$k(d')\!<\!k(c)$%
, il y a 
$e\in ]d', n[ \subset ]d, n[$
tel que
$k(e)\!=\!q_1(e)$%
.

Si
$0\!<\!\delta_1\!<\!e-d'$
on pose
$g_{1}: I\rightarrow{\Bbb R\/}$
d\'efini par~:
$$\tilde{g}_1(t)\!=\!k(e)+\int_{e}^{t}\alpha_{\delta_1}^{e}q_1'+\beta_{\delta_1}^ek'=
k+\int_e^{t}\alpha_{\delta_1}^{e}(q_1'-k')\!=\!q_1+\int_e^{t}\beta_{\delta_1}^{e}(k'-q_1')\!$$

Si
$t\!\geq\!e$%
, comme
$\alpha_{\delta_{1}}(t)\!=\!0$%
, on 
$\tilde{g}_{1}(t)\!=\!k(t)$%
, donc
$k\!\leq_{[e, n']}\tilde{g}_{1}$%
. Par convexit\'e on a
$k\!\leq_{[b', d']} \tilde{g}_{1}$
puisque
$q_{1}\!=\!k(e)+\int_{e}q_{1}',\ k\!=\!k(e)+\int_{e}k'$
et
$k\!\leq_{[b', d']} q_{1}$%
. On a aussi
$k\!\leq_{[d', e]} \tilde{g}_{1}$
puisque
$q'_{1}\!\leq_{[d', n']} k'$%
. Ainsi
$k\!\leq_{[b', n']} \tilde{g}_{1}$
et
$(2)$
sauf l'\'egalit\'e pour
$t\!\leq\!b$.
.

Si
$0\!<\!\delta\!<\!b-c'$
il suffit de prendre~:
$$g_1=\alpha_{\delta, b}k+\beta_{\delta, b}\tilde{g}_1$$
\hfill\findem
\vskip-2mm

\Demd u Lemme D|
La {\petcap Scolie\/} avec
$(k, q)\!=\!(k, p_0)$
donne une fonction
$g_1\!=\!k_1\!\geq\!k$
 en
$c, d, n$%
-dromadaire,  la {\petcap Scolie\/} pour
$(k, q)\!=\!(k_1, p_1)$
donne alors une fonction non singuli\`ere
$k_2\!\geq\!k_1\!\geq\!k$
et il suffit de poser
$$k_s\!=\!\alpha_1^1(s)[\alpha_0^1(s)k+\beta_0^1(s)k_1]+\beta_1^1(s)k_2$$
\hfill\findem
\npg
\fancyhead[RE]{Index des notations}
  \fancyhead[LE]{\thepage}
%\input DIN.tex
%\titrecourant={\hfill\\small Index des notations\hfill}
%\parc
\vskip5mm
\centerline{\petcap Indexes}

\parc
Dans la colonne de droite
${\bf 0\/}\, {\romannumeral 1}\  2$
(resp.
${\bf N\/}\ \ 2$%
)
et\ \
${^1}\ 2$
signifient \hfill\break
{\bf D\'efinition 0\/}
${\romannumeral 1})$
page 2.\ \
(resp. {\bf Notations\/}\ \
page 2.)
et\ \
note%
$^1$ page 2.
\finc
\null\vskip5mm
{\small
\centerline{\bf Index des notations}
\vskip2mm

\vskip2mm
{%
$\underline{\cal V\/}_W=({\cal V\/}, T, \pi, q)    $%\/}
\/}
\tpage {${\bf 0\/}\, {\romannumeral 1})\ 2$ }
{%
$r_W=r_{f, W}    $%\/}
\/}
\tpage {${\bf 0\/}\, {\romannumeral 1})^{2}\ 2$ }
{%
$h_W=h_{f, W}    $%\/}
\/}
\tpage {${\bf 0\/}\, {\romannumeral 1})\ 2$ }
{%
$h\leq_Y, h<_Y    $%\/}
\/}
\tpage {$\ 2$ }
{%
${_{\kappa, f}}F={_{\kappa}}F=$%\/}
\/}
\tpage {${\bf 2\/}\, \ \ 3$ }
{%
 $\bigl({\cal P\/}{_{u'}}W\bigr)
$%\/}
\/}
\tpage {$\ 3$ }
{%
${\goth A\/}(V)_v={\goth A\/}_v,\ {\goth M\/}(V)_v={\goth M\/}_v,\ {\goth A\/}(V)_v={\goth A\/}_v, 
$%\/}
\/}
\tpage {$\ 3$}
{%
${\cal X\/}(V)\!=\!{\cal X\/}, {\cal Y\/}(V)\!=\!{\cal Y\/}, {\cal X\/}_v, {\cal Y\/}_v
$%\/}
\/}
\tpage {$ \ 3$}
{%
$T_{v}V
$%\/}
\/}
\tpage {$\ 3$}
{%
${\goth J\/}_w(f)
$%\/}
\/}
\tpage {$\ 4$}
{%
${\mathop{\rm \Sigma\/}\nolimits}(f), {\mathop{\rm \Sigma\/}\nolimits}_W(f)
$%\/}
\/}
\tpage {${\bf 3\/}\, {\romannumeral 2}1)\ 4$ }
{%
$q_{H,w}
$%\/}
\/}
\tpage {${^{14}} \ 4$ }
 {%
 $\bigl({\cal D\/}\bigr)
$%\/}
\/}
\tpage {${\bf 3\/}\, {\romannumeral 1}2)\ 4  $}
{%
 $\bigl(D_{\goth W\/}\bigr), \bigl(P\bigr)
$%\/}
\/}
\tpage {${\bf 3\/}\, {\romannumeral 2}2)\ 4  $}
{%
${\mathop{\rm Ind\/}\nolimits}_{c}(f)
$%\/}
\/}
\tpage {${\bf 4\/}\, {\romannumeral 1})\ 5  $}
{%
$\mu=\mu^N\oplus\mu^P
$%\/}
\/}
\tpage {${\bf 4\/}\, {\romannumeral 4})\ 5  $}
{%
${\cal F\/}=(f, Z), {\cal F\/}=(f, \Psi,  Z)\build{=}_{}^{\rm abr.}(f, Z)
$%\/}
\/}
\tpage {${\bf 4\/}\, {\romannumeral 6})\ 5$ }
{%
$\leq^Z
$%\/}
\/}
\tpage {${\bf 5\/}\, {\romannumeral 1})\ 5$ }
{%
$\preceq^Z
$%\/}
\/}
\tpage {${\bf 5\/}\, {\romannumeral 2})\ 6$ }
{%
${c_{_\searrow}}\!(Z)$%\/}
\/}
\tpage {${\bf 5\/}\, {\romannumeral 3})\ 6$ }
{%
${^{^\nwarrow}c}(Z)$%\/}
\/}
\tpage {${\bf 5\/}\, {\romannumeral 4})\ 6$ }
{%
$\overline{c}_{_{\searrow\!\cdots\!\searrow}}$%\/}
\/}
\tpage {${\bf 5\/}\, {\romannumeral 5})\ 6$ }
{%
$\bigl(\overline{\cal P\/}{_{u'}}c\bigr)$%\/}
\/}
\tpage {${\bf \/}\, \ 6$ }
{%
$\Phi_{t}^{U}, \Lambda_{t}$%\/}
\/}
\tpage {${\bf \/}\, \ 6$ }
{%
$\theta, \chi, \chi(n, \tau)\!=\!n+\tau p +{{\tau}\over{s}}\cdot\theta({{\tau}\over{s}}\cdot n)$%\/}
\/}
\tpage {${\bf \/}\, \ 7$ }
{%
$\beta, \beta_{\epsilon, a}, \alpha_{\epsilon}^{a}, \beta_{\epsilon}^{a}, \alpha_{\epsilon, a}$%
\/}
\tpage {$ \ 18$ }
}
\npage
\npage
%\input DIT.tex
%\titrecourant={\hfill\\small Index terminologique\hfill}
%\parc
\fancyhead[RE]{Table des mati\`eres}
 \fancyhead[LE]{\thepage}
\null\vskip5mm

%{\small

{\small
\centerline{\bf Index terminologique}
\vskip3mm
{%
{%
{\it adapt\'ee {\rm [\/}m\'etrique adapt\'ee \`a {\rm ]\/}%\/}
\/}
\tpage {${\bf 4\/}\, {\romannumeral 4})\ 5$ }
{%
{\it avant {\rm [\/}point $Z$-avant un autre {\rm ]\/}
\/}
\tpage {${\bf 5\/}\, {\romannumeral 1})\ 5$ }
{\it chemin de fonctions de \`a\/}
\/}
\tpage {${\bf 1\/}\, {\romannumeral 2})\ 2$ }
{%
{\it coindice  d'un point critique\/}
\/}
\tpage {${\bf 4\/}\, {\romannumeral 1})\ 5$ }
{%
{\it conappe d'un point critique%\/}
\/}
\tpage {${\bf 5\/}\, {\romannumeral 4})\ 6$ }
{%
{\it  cr\^ete d'un val%\/}
\/}
\tpage {${\bf 0\/}\, {\romannumeral 1})\ 2$ }
{%
{\it dromadaire {\rm [\/}fonction en $c, d, n$-dromadaire {\rm ]\/}%\/}
\/}
\tpage {${\bf D\/}\, \ 24$ }
{%
{\it  ensemble de points critiques d'une fonction\/}
\/}
\tpage {${\bf 3\/}\, {\romannumeral 2}1)\ 4$ }
{%
{\it donn\'ee (de Morse)\/}
\/}
\tpage {${\bf 4\/}\, {\romannumeral 6})\ 5$ }
{%
{\it fonction de Morse\/}
\/}
\tpage {${\bf4\/}\, {\romannumeral 3})\ 5$ }
{%
{\it  fonction val {\rm ou\/} de val\/}
\/}
\tpage {${\bf 0\/}\, {\romannumeral 1})\ 2$ }
{%
{\it  fond d'un val\/}
\/}
\tpage {${\bf 0\/}\, {\romannumeral 1})\ 2$ }
{%
{\it indice  d'un point critique\/}
\/}
\tpage {${\bf 4\/}\, {\romannumeral 1})\ 5$ }
{%
{\it  infinit\'esimalment stabilis\'ee d'une fonction%\/}
\/}
\tpage {${\bf 3\/}\, {\romannumeral 3}2)\ 4$ }
{%
{\it m\'etrique adapt\'ee\/}
\/}
\tpage {${\bf 4\/}\, {\romannumeral 4})\ 5$ }
{%
{\it $Z$-modifi\'ee {\rm [\/}donn\'ee de Morse (faible) $Z$-modifi\'ee {\rm ]\/}%\/}
\/}
\tpage {${\bf 4\/}\, {\romannumeral 7})\ 5$ }
{%
{\it nappe (d'un ferm\'e critique \`a un niveau)\/}
\/}
\tpage {${\bf 2\/}\ 3$ }
{%
{\it nappe d'un point critique\/}
\/}
\tpage {${\bf 5\/}\, {\romannumeral 3})\ 6$ }
{%
{\it nappe satur\'ee d'un point critique\/}
\/}
\tpage {${\bf 5\/}\, {\romannumeral 5})\ 6$ }
{%
{\it param\'etrisation de Morse d'un point critique\/}
\/}
\tpage {${\bf 4\/}\, {\romannumeral 2})\ 5$ }
{%
{\it plat [champ de vecters
$f$%
-plat]\/}
\/}
\tpage {${\bf 3\/}\, {\romannumeral 1}0)\ 4$ }
{%
{\it plat [tube
$f$%
-plat]\/}
\/}
\tpage {${\bf 3\/}\, {\romannumeral 1}1)\ 4$ }
{%
{\it propre {\rm [\/}application $F$-propre {\rm ]\/}%\/}
\/}
\tpage {${\bf 2\/}\, \ 3$ }
{%
{\it pseudogradient (faible)\/}
\/}
\tpage {${\bf 4\/}\, {\romannumeral 5})\ 5$ }
{%
{\it rahla d'une fonction en $c, d, n$-dromadaire\/}
\/}
\tpage {${\bf D\/}\, \ 24$ }
{%
{\it sous {\rm [\/}point $Z$-sous un autre {\rm ]\/}%\/}
\/}
\tpage {${\bf 5\/}\, {\romannumeral 2})\ 6$ }
{%
{\it  stabilis\'ee d'une fonction%\/}
\/}
\tpage {${\bf 3\/}\, {\romannumeral 3}1)\ 4$ }
{%
{\it  taille d'un val%\/}
\/}
\tpage {${\bf 0\/}\, {\romannumeral 1}2)\ 2$ }
{%
{\it val de Morse\/}%\/}
\/}
\tpage {${\bf 4\/}\, {\romannumeral 9})\ 5$ }
{%
{\it val constant%\/}
\/}
\tpage {${\bf 0\/}\, {\romannumeral 1}4)\ 2$ }
{\it val pour une fonction%\/}
\/}
\tpage {${\bf 0\/}\, {\romannumeral 1})\ 2$ }
{%
{\it (micro)voisinage tubulaire m\'etrique ouvert, ferm\'e\/}
\/}
\tpage {${^2}\ 2$ }
%\titrecourant={\hfill\\small Table des mati\`eres\hfill}
\npage
\fancyhead[RE]{Table des mati\`eres}
\fancyhead[LE]{\thepage}

\centerline{\petcap  Table des mati\`eres}
\vskip3mm
%{ Abstract.\/}
 %\tpage{1}
%{ Table des mati\`eres.\/}
 %\tpage{1}

%
\vskip3mm
\vskip3 mm
{ R\'esum\'e et Abstract.\/}
 \tpage{1}

{ Guide de lecture.\/}
 \tpage{1}
\S1 {Introduction : d\'efinitions et fil des \'enonc\'es \/}
 \tpage{2}
\S2 { R\'eduction du lemme d'\'elimination aux Lemmes 
{\bf A\/}, {\bf B\/}, {\bf C\/}, {\bf D\/}.\/}
 \tpage{9}
{\S3  Commentaires bibliographiques%
}.\/}
 \tpage {10}
\quad{\bf 1925-1969\/}    Les \'ecrits de More (et avec Huebsch)
 \tpage {10}
\quad{\bf 1960-1968\/}    Le h-cobordisme de Smale et ses expositions par Milnor et Cerf
 \tpage {12}
\quad{\bf 1956-1982\/}    ${\goth M\/}{\goth J\/}^{2}$%
, singularit\'es  et m\'ethode du chemin
par Samuel, Mother, Mather%
$,\ldots$
 \tpage {13}
\quad{\bf 2012-2014\/}  Deux \og reprises\fg de la g\'eom\'etrie de Cerf
par Laudenbach
 \tpage {14}
{  R\'ef\'erences.\/}
 \tpage {15}
\vskip5mm
\centerline{ {\bf Appendices\/}%
}
%\centerline{\eightpoint Lemmes, d\'emonstration des \'enon\c c\'es et
%entre les lignes du r\'esum\'e.%
%}
%
\vskip2mm
   Notations 
 \tpage {17}
{\bf A\/}    L'ascenseur de val  
 \tpage {17}
{\bf B\/}    Stabilisation infinit\'esimale et stabilisation 
 \tpage {18}
\quad{\bf 1\/}    La m\'ethode du chemin de Moser 
 \tpage {18}
\quad{\bf 2\/}    Un lemme
${\goth M\/}{\goth J\/}^{2}$
relatif
 \tpage {19}
\quad{\bf 3\/}    D\'emonstration du {\petcap Lemme\/} B 
 \tpage {20}
{\bf  C\/}   Canonicit\'e de l'extension d'une nappe transverse 
 \tpage {22}
{\bf  D\/}   Lemme du dromadaire
 \tpage {23}
{ Index des notations.\/}
 \tpage{25}
{ Index terminologique.\/}
 \tpage{26}

{ Table des mati\`eres.\/}
 \tpage{27}

{Quatri\`eme de couverture~: d\'eveloppement du r\'esum\'e.\/}
 \tpage{28}

 \npage
 \fancyhead[RE]{Quatri\`eme de couverture}
\fancyhead[LE]{\thepage}
\centerline{Quatri\`eme de couverture~: d\'eveloppement du r\'esum\'e}
\vskip5mm
Les r\'edacteurs du livre de Milnor sur le h-cobordisme
auquel l'arbitre fait allusion font pr\'ec\'eder la preuve du lemme d'\'elimination
des deux trois lignes~:

{\sl Remark: The proof, due to Morse [11], [32], is quite formidable. Not including the technical theorem 5.6 it occupies the following 10 pages.}

Ce livre de Milnor est consid\'er\'e comme la r\'ef\'erence oblig\'ee de r\'edaction de la th\'eorie des anses \`a l'aide des fonctions de Morse.
Cependant les \og happy few\fgf, dont le premier r\'edacteur du livre de Milnor, consid\`erent que le traitement fonctionnel de Cerf%
\footnote{\small
longtemps confidentiel d\'esormais accessible gr\^ace \`a divers programmes de num\'erisation.
} avec son \og lemme des chemins \'el\'ementaires\fg utilisant habilement l'action sur les fonctions du groupe topologique produit des diff\'eomorphismes de source et de but pour plonger sans effort les mod\`eles de d\'eformation  est bien plus pertinent et propice \`a versions \`a param\`etres.

Encore plus ais\'e est le fait que, en composant au but pr\`es d'un minimum
la fonction avec un diff\'eomorphisme de
${\Bbb R\/}$%
, on peut arbitrairement baisser sa valeur.

Le pr\'esent opuscule propose de profiter de ce
\og trivial abaissement des minimaux\fg
en consid\'erant sur une vari\'et\'e
$M$
munie d'une fonction
$f$
une sous-vari\'et\'e
$N$
telle que sur chaque disque fibre d'un voisinage tubulaire la fonction
ait un minimum au point de
$N$
correspondant, ce sont les {\it fonctions  val en
$N$%
}.

Ceci permet sans peine de prolonger, en une d\'eformation de la fonction,
une d\'eformation \og par le bas\fg de la restriction
de la fonction \`a la sous-vari\'et\'e.

Et ce \`a la fois de mani\`ere continue en fonction des donn\'ees
et avec un contr\^ole de la dynamique d'un champ de pseudo-gradient associ\'e
\`a la fonction.

La m\'ethode du chemin de Moser-Mather  caract\'erise alg\'ebriquement
ces fonctions val impliquant en particulier
le lemme de forme canonique de Morse.

On introduit enfin les m\'etriques adapt\'ees \`a une fonction de Morse.

Elles permettent
une extension lisse canonique d'un germe de sous-vari\'et\'e invariante
par le gradient de la fonction pour une telle m\'etrique
et coupant, transversalement
en un arc, une nappe ascendante d'un point critique~:

L'arc est prolong\'e au del\`a du point critique et la sous-vari\'et\'e
prolong\'ee n'a qu'un seul point critique, il est d'indice 1 et situ\'e
au point critique de d\'epart.

Apr\`es avoir donn\'e notations et  d\'efinitions formelles justifiant ce yoga, 
\'enonc\'e les trois lemmes \'evoqu\'es ci-dessus ainsi qu'un quatri\`eme
r\'esolvant explicitement le probl\`eme en dimension 1,
 on donnera en une page
une d\'emonstration du lemme d'\'elimination qui
passe \`a param\`etres et contr\^ole la dynamique de
gradient.

Les lemmes  sont  \'etablis en appendice apr\`es
que des commentaires bibliographiques r\'esument en quelques
pages la litt\'erature dont l'auteur s'est inspir\'e.

\vfill

\noindent{
Metteur en sc\`ene et secr\'etaire~: Alexis Marin Bozonat\hfill\break
\null\ courriel~: alexis.charles.marin@gmail.com}

%\vskip5mm
\noindent{
Institut Fourier, UMR 5582, Laboratoire de Math\'ematiques Universit\'e  Grenoble Alpes, CS 40700, 38058 Grenoble cedex 9, France}

\end{document}